\documentclass{article}
\usepackage{amsmath,amssymb,amsthm,amscd,latexsym,}
 
\usepackage{graphicx}

\makeatletter
\renewcommand{\mod}[1]{\allowbreak \if@display \mkern 8mu \else
\mkern 5mu\fi {\operator@font mod}\,\,#1}
\makeatother

\newcommand{\bc}{\mathbb C}

 \newcommand{\bq}{\mathbb Q}
\newcommand{\br}{\mathbb R}
\newcommand{\bz}{\mathbb Z}
\newcommand{\bff}{\mathbb F}

\newcommand{\bk}{\mathbb K}
\newcommand{\bo}{\mathbb O}

\DeclareMathOperator{\rk}{rank}
\DeclareMathOperator{\ch}{ch}

\DeclareMathOperator{\discr}{{\rm discr}\,}
\newtheorem{theorem}{Theorem}
\newtheorem{proposition}[theorem]{Proposition}
\newtheorem{definition}[theorem]{Definition}

\newtheorem{lemma}[theorem]{Lemma}

\newtheorem{conjecture}[theorem]{Conjecture}

\newcommand\F{\mathcal F}
\newcommand\Hh{\mathcal H}
\newcommand\La{\mathcal L}

\newcommand\M{\mathcal M}

\begin{document}
\title{On ground fields of arithmetic hyperbolic reflection groups}
\date{}
\author{Viacheslav V. Nikulin\footnote{This paper was written with the
financial support of EPSRC, United Kingdom (grant no. EP/D061997/1)}}
\maketitle

\begin{abstract}
Using authors's methods of 1980, 1981,
some explicit finite sets of number fields containing ground fields of
arithmetic hyperbolic reflection groups are defined, and good bounds of
their degrees (over $\bq$) are obtained. For example,
degree of the ground field of any arithmetic hyperbolic reflection group
in dimension at least 6 is bounded by 120. These results could be
important for further classification.

We also formulate a mirror symmetric conjecture to finiteness of the
number of arithmetic hyperbolic reflection groups which was established
in full generality recently. 

This paper also gives corrections to my papers \cite{Nik2} (see appendix) 
and \cite{Nik7}. 
\end{abstract}

\vskip1cm \centerline{\it Dedicated to John McKay} \vskip1cm

\section{Introduction} \label{introduction}
There are only three types of
simply-connected complete Riemannian manifolds of constant curvature:
spheres, Euclidean spaces and hyperbolic
spaces. Discrete reflection groups (generated by
reflections in hyperplanes in these spaces) were defined by H.S.M.
Coxeter. He classified these groups in spheres and
Euclidean spaces \cite{Coxeter}.

There are two types of discrete reflection groups with fundamental
domain of finite volume in hyperbolic spaces: general and
arithmetic. In this paper, we consider only arithmetic hyperbolic
reflection groups.

In \cite{Vin1}, Vinberg (1967) stated and proved
a criterion for the
arithmeticity of discrete reflection groups in hyperbolic spaces
in terms of their fundamental chambers.  In particular, he introduced
the notion of the ground field of such groups. This is a totally
real algebraic number field $\bk$ of a finite degree over $\bq$. The
arithmetic reflection group $W$ is a subgroup of finite index of the
automorphism group $O(S)$ of a hyperbolic quadratic form
(a hyperbolic lattice) $S$ over the ring of integers $\bo$ of this field.
See Sect. \ref{sec:refautforms} for the exact definition. Hyperbolic lattices
$S$ having a reflection subgroup $W\subset O(S)$ of finite index are called
{\it reflective}. One can canonically associate to $S$ the hyperbolic space
$\La(S)$ of dimension $\rk S - 1$
such that groups $W\subset O^+(S)$ act in $\La(S)$ and define the
arithmetic hyperbolic reflection group $W$. One can always take the maximal
reflection subgroup $W(S)\subset O(S)$, the reflection group of the
lattice $S$,  which contains $W$. Thus, $S$ is reflective if and only if
$[O(S):W(S)]<\infty$. Two
hyperbolic lattices which differ by multiplication of their forms by
$k\in \bk$ are called {\it similar}.
Their hyperbolic spaces and automorphism groups are clearly identified.

Thus, {\it classification of similarity classes of reflective hyperbolic
lattices} is the key problem in classification of arithmetic
hyperbolic reflection groups. It includes classification
of maximal arithmetic hyperbolic reflection groups, and it is important
for classification of hyperbolic Lie algebras.
The degree $N=[\bk:\bq]$ of the
ground field $\bk$ of $S$ and $W$, and the dimension $n=\dim \La(S)=
\rk S-1\ge 2$ are the most important parameters for this classification.

In \cite{Nik1, Nik2}, the author proved (1980, 1981) that the
number of similarity classes of reflective hyperbolic lattices is
finite for the fixed parameters $n$ and $N$. Moreover, the number
of ground fields of the fixed degree $N$ is also finite. In
\cite{Nik2}, the author proved (1981) that there exists a
effective constant $N_0$ such that $N\le N_0$ if the dimension
$n\ge 10$.

In \cite{Vin2}, \cite{Vin3},  Vinberg proved (1981) that
$n < 30$: arithmetic hyperbolic reflection groups don't exist in dimensions
$n\ge 30$.

Thus, the numbers of similarity classes of reflective hyperbolic lattices
and maximal arithmetic hyperbolic reflection groups are finite in dimensions
$n\ge 10$.

About these results, see also reports \cite{Vin4} and \cite{Nik3}
at International Congresses of Mathematicians.

Almost 25 years boundedness of degree $N=[\bk:\bq]$ remained opened in
small dimensions $2\le n\le 9$. Only in 2005 it was proved in
dimension $n=2$ by Long, Maclachlan and Reid \cite{LMR} and in
dimension $n=3$ by Agol \cite{Agol}. In 2006, the author shown
\cite{Nik5} that the boundedness in remaining dimensions $4\le n \le 9$
can be easily deduced from boundedness in dimensions $n=2,\,3$ and methods
of \cite{Nik1} and \cite{Nik2} (see Theorem \ref{thfor4} and its proof
in Sect. \ref{subsec:ge4}).

Thus, now, finiteness of the numbers of reflective hyperbolic lattices and
of maximal arithmetic hyperbolic reflection groups are established in
full generality: for all dimensions of hyperbolic spaces and for all
ground fields together.

Unfortunately, these finiteness results are very far from classification of
these finite sets. The purpose of this paper, is to prove some
explicit results in this direction.
Perhaps, the first and the most important
problem is to enumerate possible ground fields $\bk$ and their degrees
$[\bk:\bq]$.

First explicit results in this direction were obtained by Vinberg
\cite{Vin2}, \cite{Vin3} in 1981. He had shown: for dimensions $n\ge 30$,
the set of ground fields is empty;
for dimensions $n\ge 22$, ground fields belong to
union of
$\F L^4=\{\bq,\,\bq(\sqrt{2}),\, \bq(\sqrt{5})\}$ and $\{ \bq(\cos(2\pi/7))\}$;
for dimensions $n\ge 14$ ground fields belong to the set $\F T$
of 13 fields which are ground fields of arithmetic triangle (plane)
reflection  groups classified by Takeuchi \cite{Tak2}, \cite{Tak3} in 1977
(their degree is bounded by 5). See Sect. \ref{subsectriangle} about
arithmetic triangle groups.

Our results can be considered as some extensions of these explicit
Vinberg's results to smaller dimensions.

In Sect. \ref{subsecedggraphs4}, we introduce finite sets of
fields $\F \Gamma^{(4)}_i(14)$, $i=1,2,3,4,5$. They are ground
fields of some V-arithmetic 3-dimensional fundamental edge chambers
of minimality 14 described by their connected hyperbolic
Gram graphs $\Gamma^{(4)}_i$, $i=1,2,3,4,5$, with four vertices
in Figures \ref{figfivegraphs} --- \ref{graphg45}.
They give some special types of V-arithmetic edge chambers introduced
and used in \cite{Nik2}. Using methods of \cite{Nik2}, in Theorem
\ref{thdegrees} and Sect. \ref{secfieldedggraphs4} we show
that the degrees of fields from these sets are bounded by
reasonable constants $24$, $39$, $53$, $120$ and $120$ respectively.

Following methods of \cite{Nik1}, \cite{Nik2} and \cite{Vin3},
in Theorem \ref{thfor10} we show that in dimensions $n\ge 10$,
the ground field of any arithmetic hyperbolic reflection group belongs
to one of sets:
$\F L^4$, $\F T$, $\F \Gamma^{(4)}_i(14)$, $i=1,2,3,4$. In particular,
the degree is bounded by $120$. Thus, author's result in \cite{Nik2} becomes
very explicit: the constant $N_0$ which we mentioned above can
be taken to be $N_0=120$.

In Sect. \ref{subsec:ge6}, we introduce one more set
$\F_{2,4}(14)$ of fields.
It is the set of ground fields of arithmetic quadrangles of minimality $14$
(i. e. ground fields of 2-dimensional arithmetic hyperbolic reflection
groups with quadrangle fundamental polygon of minimality $14$). According
to Takeuchi \cite{Tak4}, the set of ground fields of arithmetic
quadrangles is finite, and degrees of these fields are bounded by $11$.

Following methods of \cite{Nik1}, \cite{Nik2}, in Theorem \ref{thfor6}
we show that in dimensions $n\ge 6$
the ground field of any arithmetic hyperbolic reflection group
belongs to one of sets
$\F L^4$, $\F T$, $\F \Gamma^{(4)}_i(14)$, $i=1,2,3,4,5$
and $\F_{2,4}(14)$.
In particular, its degree is bounded by 120.

Unfortunately, now we don't have so explicit results for smaller dimension.
Following \cite{Nik5}, we only show in Theorem
\ref{thfor4} that in dimensions $n\ge 4$, the ground
field of any arithmetic hyperbolic reflection group belongs to one of
sets $\F L^4$, $\F T$, $\F \Gamma^{(4)}_i(14)$, $i=1,2,3,4$, or it is
the ground field of $3$ or $2$-dimensional arithmetic hyperbolic
reflection group with a fundamental chamber of minimality $14$. Thus, the
degree is bounded by the maximum of 120 and of degrees of ground fields of
$3$ and $2$-dimensional arithmetic hyperbolic reflection groups of
minimality $14$.

For 3-dimensional arithmetic hyperbolic reflection groups, we know finiteness
of the number of ground fields by Agol \cite{Agol}. By author's knowledge,
no explicit bound of their degree is known.

Following Long, Maclachlan and Reid \cite{LMR} and
Borel \cite{Bor}, Takeuchi \cite{Tak4},
in Sect. \ref{subsec:deg=2} we show that degree of the ground field
of any 2-dimensional arithmetic hyperbolic reflection group is bounded by
$44$.

\medskip

Using known at that time finiteness result for hyperbolic reflective lattices
over $\bz$, in \cite{GrNik1}--\cite{Nik4} some finiteness results for
IV type (i. e. of signature $(2,t)$) integer reflective lattices $S$ were
obtained, and some general conjectures were formulated. Here a lattice
$S$ over $\bz$ of signature $(2,t)$ is called reflective if IV type Hermitian
symmetric domain associated to $S$ has an $O^+(S)$-automorphic
holomorphic form $\Phi$ of positive weight such that all components
of its divisor are  quadratic divisors orthogonal to roots
(giving reflections) of $S$.
This automorphic form is called {\it reflective.}
One can consider these finiteness statements about IV type reflective
lattices over $\bz$
as ``mirror symmetric'' to finiteness results about hyperbolic reflective
lattices over $\bz$. Since we now know finiteness of hyperbolic reflective
lattices in general, in  Sect. \ref{sec:refautforms}
we formulate the corresponding mirror symmetric conjecture about IV type
reflective lattices in general --- over arbitrary totally real algebraic
number fields. We expect that it is valid.

Some arithmetic hyperbolic reflection groups and some reflective automorphic
forms and corresponding hyperbolic and IV type reflective lattices
over $\bz$ are important in Borcherds proof \cite{Bor}
of Moonshine Conjecture by Conway and Norton \cite{CN}
which had been first discovered by John McKay.

We hope that similar objects over arbitrary number fields will find similar
astonishing applications in the future. At least, the results and conjectures
of this paper show that they are very exceptional even in this very
general setting.

\medskip

At first, the paper appeared as preprint \cite{Nik6} which   
was published in \cite{Nik7}. In Appendix 
we review and correct Section 1 of our old paper \cite{Nik2} 
which was used in these papers. The present variant takes these 
corrections under considerations.  

\section{Reminding of some basic facts about hyperbolic fundamental
polyhedra}\label{sec:V-arpol}

Here we remind some basic definitions and results about fundamental
chambers (always for discrete reflection groups) in hyperbolic spaces and
their Gram matrices.
See \cite{Vin1}, \cite{Vin5} and \cite{Nik1}, \cite{Nik2}.

We work with Klein model of a hyperbolic space
$\La$ associated to a hyperbolic form $\Phi$ over the field of
real numbers $\br$ with signature $(1,n)$, where $n=\dim \La$.
Let $V=\{x\in \Phi|x^2>0\}$ be the cone determined by $\Phi$,
and let $V^+$ be one of the two halves of this cone. Then
$\La=\La(\Phi)=V^+/\br^+$ is the set of rays in $V^+$; we let
$[x]$ denote the element of $\La$ determined by the ray $\br^+x$
where $x\in V^+$ and $\br^+$ is the set of all positive real numbers.
The hyperbolic distance is given by the formula
$$
\rho([x],[y])=(x\cdot y)/\sqrt{x^2y^2},\ [x], [y]\in \La,
$$
then the curvature of $\La$ is equal to $-1$.

Every half-space $\Hh^+$ in $\La$
determines and is determined by the orthogonal element $e\in \Phi$ with
square $e^2=-2$:
$$
\Hh^+=\Hh_e^+=\{[x]\in \La|x\cdot e\ge 0\}.
$$
It is bounded by the hyperplane
$$
\Hh^+=\Hh_e^+=\{[x]\in \La|x\cdot e = 0\}
$$
orthogonal to $e$.  If two half-spaces
$\Hh_{e_1}^+$, $\Hh_{e_2}^+$ where $e_1^2=e_2^2=-2$
have a common non-empty open subset in $\La$, then
$\Hh_{e_1}\cap \Hh_{e_2}$ is an angle of the value $\phi$ where
$2\cos{\phi}=e_1\cdot e_2$ if $-2<e_1\cdot e_2\le 2$,
and the distance between hyperplanes $\Hh_{e_1}$ and
$\Hh_{e_2}$ is equal to $\rho$ where $2\ch{\rho}=e_1\cdot e_2$
if $e_1\cdot e_2>2$.

A convex polyhedron $\M$ in $\La$ is intersection of a finite
number of half-spaces $\Hh^+_e$, $e\in P(\M)$, where $P(\M)$ are all the
vectors with square $-2$ which are orthogonal to the faces
(of the codimension one) of $\M$ and are directed outward. The matrix
\begin{equation}
A=(a_{ij})=(e_i\cdot e_j),\ e_i, e_j\in P(\M),
\label{Grammatr}
\end{equation}
is the Gram matrix $\Gamma (\M)=\Gamma (P(\M))$ of $\M$.
It determines $\M$ uniquely up to motions
of $\La$. If $\M$ is sufficiently general, then $P(\M)$ generates $\Phi$,
and the form $\Phi$ is
\begin{equation}
\Phi=\sum_{e_i,e_j\in P(\M)}{a_{ij}X_iY_j}\mod {\rm Kernel},
\label{formPhi}
\end{equation}
and $P(\M)$ naturally identifies with a subset of $\Phi$ and defines $\M$.

The polyhedron $\M$ is a fundamental chamber of a discrete reflection group
$W$ in $\La$ if and only if $a_{ij}\ge 0$ and
$a_{ij}=2\cos {\frac{\pi}{m_{ij}}}$ where $m_{ij}\ge 2$ is an integer if
$a_{ij}<2$ for all $i\not=j$.
Symmetric real matrices $A$ satisfying these conditions and having all their
diagonal elements equal to $-2$ are called {\it fundamental} (then the set
$P(\M)$ formally corresponds to indices of the matrix $A$).
As usual, further we identify  fundamental matrices with fundamental
graphs $\Gamma$. Their vertices correspond to
$P(\M)$. Two different vertices $e_i\not=e_j\in P(\M)$ are connected by
the thin edge of the integer weight $m_{ij}\ge 3$ if
$0<a_{ij}=2\cos {\frac{\pi}{m_{ij}}}<2$, by the thick edge if $a_{ij}=2$,
and by the broken edge of the weight $a_{ij}$ if $a_{ij}>2$.
In particular, the vertices $e_i$ and $e_j$ are disjoint if and only if
$e_i\cdot e_j=a_{ij}=2\cos{\frac{\pi}{2}}=0$. Equivalently, $e_i$ and
$e_j$ are perpendicular (or orthogonal). See some examples of such graphs
in Figures \ref{figLanner} --- \ref{graphg45} below.

For a real $t>0$, we say that a fundamental matrix $A=(a_{ij})$ (and the
corresponding fundamental chamber $\M$) {\it has minimality $t$}
if $a_{ij}<t$ for all $a_{ij}$. Here we follow \cite{Nik1},
\cite{Nik2}. Further, the minimality $t=14$ will be especially important.

It is known that fundamental domains of arithmetic hyperbolic
groups must have finite volume. Let us assume that it is valid for a
fundamental chamber $\M$ of a hyperbolic discrete reflection group.
As Vinberg had shown \cite{Vin1}, in order for $\M$ to be a fundamental
chamber of an arithmetic reflection group $W$ in $\La$, it is necessary and
sufficient that
all of the cyclic products
\begin{equation}
b_{i_1\dots i_m}=a_{i_1i_2}\cdot a_{i_2i_3}\cdots a_{i_{m-1}i_m}\cdot
a_{i_mi_1}
\label{cycprod}
\end{equation}
be algebraic integers, that the field
$\widetilde{\bk}=\bq(\{a_{ij}\})$ be
totally real, and that, for any embedding $\widetilde{\bk}\to \br$ not
the identity over the {\it ground field} $\bk=\bq(\{b_{i_1\dots i_m}\})$
generated by all of the cyclic products \eqref{cycprod}, the form
\eqref{formPhi}
be negative definite.

Fundamental real matrices $A=(a_{ij})$, $a_{ij}=e_i\cdot e_j$,
$e_i,\,e_j\in P(\M)$ (or the corresponding graphs),
with a hyperbolic form $\Phi$ in \eqref{formPhi} and satisfying
these Vinberg's conditions will be further called
{\it V-arithmetic} (here we don't require that the corresponding
hyperbolic polyhedron $\M$ has finite volume). It is well-known
(and easy to see; see arguments in Sect. \ref{subsec:gam(4)1}) that
a subset $P\subset P(\M)$ also defines a V-arithmetic
matrix $(e_i\cdot e_j)$, $e_i, e_j\in P$, with the same ground
field $\bk$ if the subset $P$ is hyperbolic, i. e. the corresponding to
$P$ form \eqref{formPhi} is hyperbolic.

\section{V-arithmetic edge polyhedra}{\label{V-arithedgepol}

A fundamental chamber $\M$ (and the corresponding Gram matrix $A$ or
a graph) is called {\it edge chamber (matrix, graph)} if all hyperplanes
$\Hh_e$, $e\in P(\M)$, contain one of two distinct vertices $v_1$ and
$v_2$ of the 1-dimensional edge $v_1v_2$ of $\M$. Assume that both vertices
$v_1$ and $v_2$ are finite (further we always consider this case).
Further we call this edge chambers {\it finite}. Assume that
$\dim \La=n$. Then $P(\M)$ consists of
$n+1$ elements: $e_1$, $e_2$ and $n-1$ elements $P(\M)-\{e_1,e_2\}$.
Here $P(\M)-\{e_1,e_2\}$ corresponds to hyperplanes which contain the edge
$v_1v_2$ of $\M$. The $e_1$ corresponds to the hyperplane which contains
$v_1$ and does not contain $v_2$. The $e_2$ corresponds to the hyperplane
which contains $v_2$ and does not contain $v_1$. Then the set $P(\M)$
is hyperbolic (it has hyperbolic Gram matrix),
but its subsets $P(\M)-\{e_1\}$ and $P(\M)-\{e_2\}$ are negative definite
(they have negative definite Gram matrix) and define Coxeter graphs.
Only the element $u=e_1\cdot e_2$ of the Gram matrix of $\M$ can be greater
than $2$. Thus, $\M$ will have the minimality $t$ if and only if
$u=e_1\cdot e_2<t$.

From considerations above, the Gram graph $\Gamma (P(\M))$ of an edge chamber
has only one hyperbolic connected component $P(\M)^{hyp}$ (containing
$e_1$ and $e_2$) and several negative definite connected components.
Gram matrix $\Gamma (P(\M)^{hyp})$ evidently also corresponds to an
edge chamber of the dimension $\#P(\M)^{hyp}-1$. If  $\M$ is
V-arithmetic, the ground field $\bk$ of $\M$ is the same as for
the hyperbolic connected component $\Gamma (P(\M)^{hyp})$

The following result had been proved in \cite{Nik2}.

\begin{theorem} (\cite[Theorem 2.3.1]{Nik2}) Given $t>0$, there exists an
effective constant $N(t)$ such that every
V-arithmetic edge chamber of the minimality $t$ with
ground field $\bk$ of degree greater that $N(t)$ over $\bq$ has the hyperbolic
connected component of its Gram graph
which has less than $4$ elements.
\label{tholdVedge}
\end{theorem}

Considerations in \cite{Nik2} (and also \cite{Nik1})
also show that the set of possible
ground fields $\bk$ of hyperbolic
connected components with at least $4$ vertices
of V-arithmetic edge chambers of minimality $t$ is finite.
Even the set of Gram graphs
$\Gamma(P(\M)^{hyp})$ of minimality $t$
with fixed $\ge 4$ number of vertices is finite.
Taking this under consideration, here we want to formulate and prove
more efficient variant of this theorem. We restrict by the minimality
$t=14$ to get an exact estimate for the constant $N(14)$, but the
same finiteness results are valid for any $t>0$.

To formulate this new variant, we need to introduce some fundamental
edge graphs.

\begin{figure}
\begin{center}
\includegraphics[width=6cm]{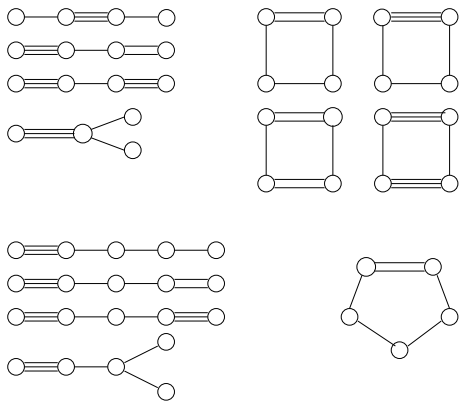}
\end{center}
\caption{All arithmetic Lann\'er graphs with at least $4$ vertices}
\label{figLanner}
\end{figure}

\subsection{Arithmetic Lann\'er graphs with $\ge 4$ elements}
\label{subsecLanner}

We remind that Lann\'er graphs are Gram graphs of bounded
fundamental hyperbolic simplexes. They are characterized as
hyperbolic fundamental graphs such that any their proper subgraph
is a Coxeter graph. They were classified by Lann\'er \cite{Lan}.
In Figure \ref{figLanner} we give all arithmetic Lann\'er graphs
with at least $4$ vertices (only one Lann\'er graph with $\ge 4$
vertices is not arithmetic). As usual, we replace a thin edge of
the weight $k$ by $k-2$-edges for a small $k$. Ground fields of
Lann\'er graphs with $\ge 4$ vertices give three fields:
\begin{equation}
\F L^4=\{\bq,\,\bq(\sqrt{2}),\ \bq(\sqrt{5})\}.
\label{fLanner4}
\end{equation}
See \cite{Vin3} for details.

\subsection{Arithmetic triangle graphs}
\label{subsectriangle}

\begin{figure}
\begin{center}
\includegraphics[width=6cm]{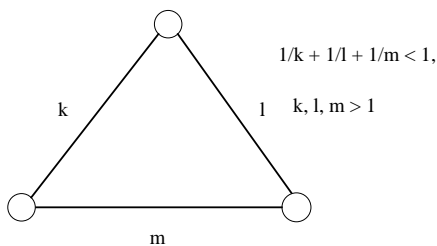}
\end{center}
\caption{Triangle graphs}
\label{figtriangles}
\end{figure}

Triangle graphs are Gram graphs of bounded fundamental triangles on
hyperbolic plane (we don't consider non-bounded triangles).
Equivalently, they are Lann\'er graphs with $3$
vertices. They are given in Figure \ref{figtriangles} where
$2\le k,l,m$ and
$$
\frac{1}{k}+\frac{1}{l}+\frac{1}{m}<1.
$$
Arithmetic triangles were enumerated
by Takeuchi \cite{Tak2}. All bounded arithmetic triangles
are given by the following triplets $(k,l,m)$:
$$
(2,3,7\ -\ 12),\ (2,3,14),\ (2,3,16),\ (2,3,18),\ (2,3,24),\ (2,3,30),\
(2,4,5\ -\ 8),
$$
$$
(2,4,10),\ (2,4,12),\ (2,4,18),\ (2,5,5),\ (2,5,6),\ (2,5,8),\
(2,5,10),\ (2,5,20),
$$
$$
(2,5,30),\ (2,6,6),\ (2,6,8),\ (2,6,12),\ (2,7,7),\ (2,7,14),\
(2,8,8),\
(2,8,16),
$$
$$
(2,9,18),\ (2,10,10),\ (2,12,12),\ (2,12,24),\ (2,15,30),\ (2,18,18),\
(3,3,4\ -\ 9),\
$$
$$
(3,3,12),\ (3,3,15),\ (3,4,4),\ (3,4,6),\ (3,4,12),\ (3,5,5),\ (3,6,6),\
(3,6,18),
$$
$$
(3,8,8),\ (3,8,24),\ (3,10,30),\ (3,12,12),\ (4,4,4\ - 6),\ (4,4,9),\
(4,5,5),\ (4,6,6),
$$
$$
(4,8,8),\ (4,16,16),\ (5,5,5),\ (5,5,10),\ (5,5,15),\ (5,10,10),\
(6,6,6),\ (6,12,12),
$$
$$
(6,24,24),\ (7,7,7),\ (8,8,8),\ (9,9,9),\ (9,18,18),\ (12,12,12),\ (15,15,15).
$$
Their ground fields were found by Takeuchi \cite{Tak3}. They give the set
of fields
\begin{equation}
\begin{array}{l}
\F T=
\{\bq\}\cup \{\bq(\sqrt{a})\ |\ a=2,\,3,\,5,\,6\}\cup
\{\bq(\sqrt{2},\sqrt{3}),\ \bq(\sqrt{2},\sqrt{5})\}\cup\\
\cup \{\bq(\cos{\frac{2\pi}{b}})\ |\ b=7,\,9,\,11,\,15,\,16,\,20\}.
\end{array}
\label{ftriangles}
\end{equation}

\subsection{V-arithmetic connected finite edge graphs with $4$ vertices
for $2<u<14$}\label{subsecedggraphs4}

Using classification of Coxeter graphs, it is easy to draw
all possible pictures of connected finite edge graphs $\Gamma^{(4)}$
with $4$ vertices and $u=e_1\cdot e_2>2$. They correspond to all
3-dimensional finite fundamental edge polyhedra with connected Gram
graph and $u>2$. They are given in Figure \ref{figfivegraphs} and give
five types of graphs $\Gamma=\Gamma^{(4)}_i$, $i=1,\,2,\,3,\,4,\,5$.
All possible natural parameters $s,\,k,\,r,\,p\ge 2$ for these graphs can
be easily enumerated by the condition that
$\Gamma-\{ e_1 \}$, $\Gamma-\{ e_2\}$ are Coxeter graphs.
They will be given in Sec. \ref{secfieldedggraphs4}  below.

\begin{definition}
For $i=1,\,2,\,3,\,4,\,5$ and $t>0$ we denote by $\Gamma^{(4)}_i(t)$ the set of
all V-arithmetic connected finite
edge graphs with 4
vertices  $\Gamma^{(4)}_i$ of the minimality $t$, i. e. for $2<u<t$, and by
$$
\F \Gamma^{(4)}_i(t)
$$
the set of all their ground fields.
\label{fieldsedge4}
\end{definition}

All V-arithmetic graphs $\Gamma^{(4)}_i$ for $2<u<t$  give particular
cases of graphs of V-arithmetic edge polyhedra
with hyperbolic connected component having $4$ vertices and minimality $t$.
Thus, by Theorem \ref{tholdVedge}, degree (over $\bq$) of fields from
$\F\Gamma^{(4)}_i(t)$ is bounded by the effective
constant $N(t)$. It follows that the sets of V-arithmetic graphs
$\Gamma^{(4)}_i(t)$ and fields $\F \Gamma^{(4)}_i(t)$ are
also finite.

Vice versa, Theorem \ref{fieldsedge4} can be deduced from finiteness
of the sets of fields above because of the
following easy statement.

\begin{figure}
\begin{center}
\includegraphics[width=10cm]{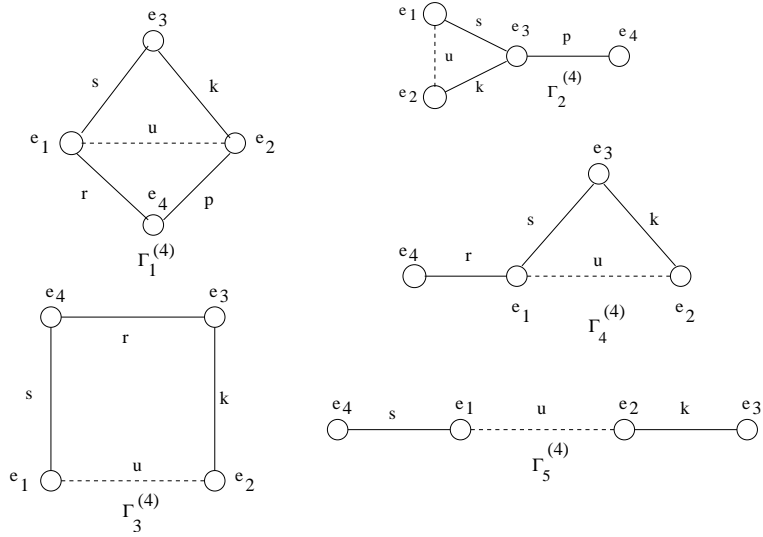}
\end{center}
\caption{Five graphs $\Gamma^{(4)}_i$, $i=1,\,2,\,3,\,4,\,5$.}
\label{figfivegraphs}
\end{figure}

\begin{proposition} The ground field of any
V-arithmetic edge chamber of the minimality $t>0$ with the hyperbolic
connected component of its Gram graph  having at least $4$ vertices
belongs to one of the finite sets of
fields $\F L^4$, $\F T$ and $\F\Gamma^{(4)}_i(t)$, $1\le i\le 5$,
introduced above.

In particular, Theorem \ref{tholdVedge} is equivalent to finiteness of
the sets of fields
$\F\Gamma^{(4)}_i(t)$, $ i=1,\,2,\,3,\,4,\,5$.
\label{prnewVedge}
\end{proposition}

\begin{proof} Let $\M$ be a V-arithmetic edge chamber of the minimality
$t>0$, and $\Gamma$ its Gram graph having the hyperbolic connected component
$\Gamma^{hyp}$ with at least 4
vertices. The graph $\Gamma^{hyp}$ is also V-arithmetic edge graph.

Assume that $u=e_1\cdot e_2<2$. Then $\{e_1, e_2\}$ gives a negative
definite subgraph of $\Gamma^{hyp}$. It
follows that any subgraph of $\Gamma^{hyp}$ is also negative definite
if it has one or two vertices. Since
$\Gamma^{hyp}$ is hyperbolic, it follows that $\Gamma^{hyp}$ contains
a minimal hyperbolic subgraph $L$ which is
Lann\'er with at least 3 vertices. Since $L$ is hyperbolic, the ground field
of $\Gamma$ is equal to the ground
field of $L$. If $L$ has more than 3 vertices, then the ground field of
$L$ is one of fields
$\F L^4$. If $L$ has 3 vertices, then the ground field of $L$ is one
of fields from $\F T$.

If $u=e_1\cdot e_2=2$, then the ground field of $\Gamma$ is equal to $\bq$.

Assume that $2<u=e_1\cdot e_2<t$. Then the subset $\{e_1, e_2\}$ is
hyperbolic and connected.
Since $\Gamma^{hyp}$ is connected, contains $e_1,\,e_2$ and
has at least 4 vertices, obviously there exists a connected subgraph
$\Gamma^{(4)}$ of $\Gamma^{hyp}$ which contains
$e_1,\,e_2$ and has four vertices.  It is hyperbolic since it contains
a hyperbolic subset $\{e_1,\,e_2\}$. Then $\Gamma^{(4)}$
is one of the hyperbolic graphs $\Gamma^{(4)}_i$, $i=1,2,3,4,5$.
Then the ground field of $\Gamma$ is equal to one of fields
$\Gamma^{(4)}_i(t)$, $1\le i\le 5$.

This finishes the proof.
\end{proof}

Degree of fields from $\F L^4$ is  bounded by $2$, and degree of fields from
$\F T$ is bounded by $5$.

For arithmetic hyperbolic reflection groups, the minimality $t=14$ is
especially important.
Using the same methods as for the proof of Theorem \ref{tholdVedge}
in \cite{Nik2},  we can prove
the following effective estimates.

\begin{theorem} The degree of fields from $\F \Gamma^{(4)}_1(14)$ is
bounded ($\le$) by $24$.

The degree of fields from $\F \Gamma^{(4)}_2(14)$ is bounded by $39$.

The degree of fields from $\F \Gamma^{(4)}_3(14)$ is bounded by $53$.

The degree of fields from $\F \Gamma^{(4)}_4(14)$ is bounded by $120$.

The degree of fields from $\F \Gamma^{(4)}_5(14)$ is bounded by $120$.

Thus, the constant $N(14)$ of Theorem \ref{tholdVedge} can be taken to
be $N(14)=120$.
\label{thdegrees}
\end{theorem}

\begin{proof} The proof requires long considerations and calculations.
It will be given later in the special Sect. \ref{secfieldedggraphs4}.
\end{proof}

In the next section, we shall consider applications of these
explicit estimates to arithmetic hyperbolic reflection groups.

\section{Application to ground fields of arithmetic hyperbolic reflection
groups}
\label{sec:ge6}

Let $W$ be an arithmetic hyperbolic
reflection group of dimension $n\ge 2$ and $\M$
its fundamental chamber. If $\M$ is not bounded, then
the ground field of $W$ is $\bq$ (see \cite{Vin1}). Since we are
interested in possible ground fields of $W$, further we
assume that $\M$ is bounded.  Then any edge $r=v_1v_2$ of $\M$
defines a V-arithmetic finite edge chamber $\M(r)$ which is intersection of all
half-spaces $\Hh_\delta^+$, $\delta\in P(\M)$, such that the
hyperplane $\Hh_\delta$ contains one of vertices $v_1$, $v_2$.
The corresponding V-arithmetic edge graph
is the Gram graph $\Gamma(r)=\Gamma(\M(r))$ of these elements
$\delta\in P(\M)$.

By \cite{Nik1} and \cite{Nik2}, there exists $e\in P(\M)$ which
defines a narrow face $\M_e=\Hh_e\cap \M$
(or a face of minimality $14$) of $\M$ (the same is valid for any hyperbolic
closed convex polyhedron). It means that for the set
$$
P(\M,e)=\{\delta\in P(\M)\ |\ \Hh_\delta\cap \Hh_e\not=\emptyset\}
$$
of neighbouring to $\M_e$ faces of $\M$ one has
$$
\delta_1\cdot \delta_2<14\ \ \forall \ \delta_1,\,\delta_2\in P(\M,e).
$$
By considering edges $r=v_1v_2$ in $\M_e$, we obtain many V-arithmetic edge
chambers $\M(r)$ and V-arithmetic edge graphs of the minimality $14$.

\subsection{Ground fields of arithmetic hyperbolic reflection \\
groups in dimensions  $n\ge 10$}\label{subsec:ge10}
Let us assume that dimension $n\ge 10$.
It was shown in  \cite{Nik2}, that  the narrow face $\M_e$ contains
an edge $r=v_1v_2$ such
that the corresponding edge chamber $\M(r)$ and its Gram
graph $\Gamma=\Gamma(r)$ (they have minimality $14$) have the hyperbolic
connected component $\Gamma^{hyp}$ with at least 4 vertices.
By Theorem \ref{tholdVedge}, then the degree of ground field of $W$ is
bounded by the constant $N(14)$.
If we additionally apply  Proposition \ref{prnewVedge}, we obtain
that the ground field of
$W$ belongs to one of finite sets of fields  $\F L^4$, $\F T$ and
$\F\Gamma^{(4)}_i(14)$, $1\le i\le 5$.
By Theorem \ref{thdegrees}, the degree of this field is bounded by $120$.

Here we want to give another proof of this result which permits to
avoid enumeration of combinatorial types
of 3-dimensional polyhedra with small number of vertices and excludes
the most difficult and large set of fields
$\F\Gamma^{(4)}_5(14)$ from the statement. This proof follows \cite{Nik4}.
It additionally uses some important arguments by Vinberg
from \cite{Vin4}.
In \cite{Vin4}, Vinberg has shown that for $n\ge 30$ there are no arithmetic
hyperbolic reflection groups; for $n\ge 22$ their ground fields
 belong to the set $\{\bq,\, \bq(\sqrt{2}), \bq(\sqrt{5}),
\bq(\cos{\frac{2\pi}{7}})\}$;
for $n\ge 14$ their ground fields belong to $\F T$. Thus,
Theorem \ref{thfor10} below
can be viewed as some extension of these statements for $n\ge 10$.

\begin{theorem} In dimensions $n\ge 10$, the ground field of any
arithmetic hyperbolic reflection group belongs to one of finite
sets of fields $\F L^4$, $\F T$, $\F\Gamma^{(4)}_1(14)$,
$\F\Gamma^{(4)}_2(14)$, $\F\Gamma^{(4)}_3(14)$ or
$\F \Gamma^{(4)}_4(14)$. In particular, its degree is
bounded ($\,\le \,$) by $120$.
\label{thfor10}
\end{theorem}

\begin{proof} Let $W$ be an arithmetic hyperbolic reflection group
of dimension $n\ge 2$ and $\M$ its fundamental chamber with
the Gram graph $\Gamma (P(\M))$. If $\M$ is not compact, then the
ground field is $\bq$. Thus, further we can assume that
$\M$ is compact. It is known, \cite{Vin1}, that then $\M$ is a
simple polyhedron which means that $\M$ is simplicial in its
vertices: any vertex is contained in exactly $n$ hyperplanes
$\Hh_\delta$, $\delta\in P(\M)$.

Arguing like in the proof of Proposition \ref{prnewVedge}, we obtain

\begin{lemma} Let $W$ be an arithmetic hyperbolic reflection group,
$\M$ its fundamental chamber and their ground
field is different from fields from $\F L^4$ and $\F T$.

Then any edge $r=v_1v_2$ of $\M$ defines an edge polyhedron $\M(r)$
such that the corresponding vertices
$e_1$ and $e_2$ of its Gram graph (that is the hyperplanes
$\Hh_{e_1}$ and $\Hh_{e_2}$ contain only
vertices $v_1$ and $v_2$ of the edge respectively)
are joined by a broken edge (i.e.  $u=e_1\cdot e_2>2$).
\label{lembroken}
\end{lemma}

For a vertex $v$ of $\M$ we denote by $C(v)$ the Coxeter graph
of $v$ which is Gram graph of all $\delta \in P(\M)$  such
that the hyperplane $\Hh_\delta$ contains $v$. Since $\M$ is
compact, $C(v)$ has exactly $n$ vertices.

\begin{lemma} Let $W$ be an arithmetic hyperbolic reflection
group, $\M$ its fundamental chamber and their ground
field is different from fields from $\F L^4$, $\F T$ and
$\F \Gamma^{(4)}_i(14)$, $i=1,\,2,\,3,\,4$.
Let $\M_e$, $e\in P(\M)$, be a narrow face of $\M$ of the minimality $14$.

Then for any vertex
$v$ of $\M$ which is contained in $\M_e$, all connected components of
the Coxeter graph $C(v)$ of $v$ have one or two vertices only.
\label{lemcoxeter}
\end{lemma}

\begin{proof}
Let $e_1$ be a vertex of $C(v)$ which is different from $e$. Let $r$
be an edge of $\M$ defined by $C(v)-\{e_1\}$. It means that all hyperplanes
$\Hh_\delta$, $\delta\in P(\M)$, which contain the edge $r$,
belong to vertices of $C(v)-\{e_1\}$. Since $e\in C(v)-\{e_1\}$, the edge $r$
belongs to $\M_e$. One of terminals of $r$ is $v$, and
$\Hh_{e_1}$ contains $v$, and it does not contain $r$. Let $v_2\in \M_e$ be
another terminal of $r$ and $e_2\in P(\M)$ gives the hyperplane
$\Hh_{e_2}$ which contains only $v_2$ and does not contain $v$. Thus,
$\Gamma(r)=C(v)\cup \{e_2\}$ is the Gram graph of the edge polyhedron of $\M$
corresponding to $r=vv_2$. By Lemma \ref{lembroken}, $e_1,e_2$ give the
only broken edge of $\Gamma(r)$.

Since $r\subset \M_e$ and $\M_e$ has minimality $14$, then $\Gamma(r)$ also
has minimality $14$.

Let us assume that the connected component of $e_1$ in
Coxeter graph $C(v)$ has more than two vertices. Then there exists a
connected subgraph $\Gamma_{e_1}$ of $C(v)$ which has three vertices and
contains $e_1$. Then $\Gamma_{e_1}\cup \{e_2\}$ is one of graphs
$\Gamma^{(4)}_i$,
$i=1,2,3,4,5$, of minimality 14.
It cannot be equal to $\Gamma^{(4)}_5$ because the connected components
of $e_1$ in $\Gamma^{(4)}_5-\{e_2\}$ and $e_2$ in $\Gamma^{(4)}_5-\{e_1\}$
have two vertices. It follows that
the ground field of $\M$ is one of fields $\F \Gamma^{(4)}_i(14)$,
$i=1,2,3,4$. We obtain a contradiction.

Thus, any vertex $e_1\not=e$ of the graph $C(v)$ has the connected
component with one or two vertices.
Obviously, the same will be valid for $e$ as well.

This finishes the proof.
\end{proof}

The crucial topological argument of the proof of Theorem \ref{thfor10}
(additional to Theorem \ref{tholdVedge} and existence of a narrow
face of minimality $14$) is as follows.
By \cite[Theorem 3.2.1]{Nik2}, for any $0\le i\le k$,
$2\le k$  and $2k-1\le m$, the average number $\alpha_m^{(i,k)}$
of $i$-dimensional faces in $k$-dimensional faces of any  $m$-dimensional
simple convex polyhedron satisfies the inequality
\begin{equation}
\alpha_m^{(i,k)}<
\frac{C_{m-i}^{k-i}\left(C^{i}_{[m/2]}+C_{m-[m/2]}^i\right)}
{C_{[m/2]}^{k}+C_{m-[m/2]}^k} .
\label{aver,i,k}
\end{equation}
In particular, for $n\ge 4$, the average number $\alpha_{n-1}^{(0,2)}$ of
vertices of $2$-dimensional faces of
a narrow face $\M_e$ (of dimension $m=n-1$) satisfies
\begin{equation}
\alpha_{n-1}^{(0,2)}<4+
\left\{
\begin{array}{rcl}
\frac{4}{n-2} & {\rm if } &\  n\ {\rm is\ even,}\\
\frac{4}{n-3} & {\rm if}   &\  n\ {\rm is\ odd.}
\end{array}
\right.
\label{alpha02}
\end{equation}

Now let us assume that the ground field is different from the fields from the
sets $\F L$, $\F T$, $\F\Gamma^{(4)}_i(14)$, $i=1,2,3,4$.
Following Vinberg \cite{Vin3},  let us estimate the number $A$ of
non-right (i.e $\not=\pi/2$) 2-dimensional angles of $\M_e$.

Let $v$ be a vertex of $\M_e$ and $C(v)$ the Coxeter graph of $v$ in $\M$.
Any 2-dimensional angle of $\M_e$ with the vertex $v$ is defined by a subset
of two distinct vertices $\{\delta_1,\delta_2\}\subset C(v)-\{e\}$. Really,
$C(v)-\{\delta_1,\delta_2\}$ define perpendicular vectors to hyperplanes
containing the plane of the angle, and
$C(v)-\{\delta_1\}$, $C(v)-\{\delta_2\}$ similarly define edges of the angle.
It is easy to see that the angle is not right if and only if $\delta_1$ and
$\delta_2$ belong to one connected component of the graph $C(v)$.
Thus, the number $A_v$ of non-right
2-dimensional angles of $\M_e$ with the vertex $v$ is equal to the number
of subsets of two distinct vertices of $C(v)-\{e\}$ which belong to one
connected component of the graph $C(v)$. By Lemma \ref{lemcoxeter},
$A_v\le [(n-1)/2]$. Thus,
\begin{equation}
\left[\frac{n-1}{2}\right]\alpha_0\ \ge\ A
\label{Abyvertices}
\end{equation}
where $\alpha_0$ is the number of vertices of $\M_e$.

Each hyperbolic triangle has at least two non-right angles. Each hyperbolic
quadrangle has at least one non-right angle.
Denoting by $\alpha_2^l$ the number of 2-dimensional faces of $\M_e$
with $l$ vertices and by
$\alpha_2$ the number of all 2-dimensional faces of $\M_e$, we obtain
\begin{equation}
A\ge \sum_{l\ge 3}{(5-l)\alpha_2^l}=5\alpha_2-\sum_{l\ge 3}{l\alpha_2^l}=
5\alpha_2-\alpha_{n-1}^{(0,2)}\alpha_2=(5-\alpha_{n-1}^{(0,2)})\alpha_2.
\label{Aby2faces}
\end{equation}

Since $\M_e$ is a simple and $(n-1)$-dimensional convex polyhedron, we have
\begin{equation}
\frac{\alpha_0(n-1)(n-2)}{2}=\alpha_2\alpha_{n-1}^{(0,2)}.
\label{alpha0alpha2}
\end{equation}

From \eqref{Abyvertices}, \eqref{Aby2faces} and \eqref{alpha0alpha2},
we obtain
$$
\alpha_{n-1}^{(0,2)}\left(\frac{(n-1)(n-2)}{2}+\left[\frac{n-1}{2}
\right]\right)\ge
\frac{5(n-1)(n-2)}{2}\ .
$$
From \eqref{alpha02}, we get
$$
\left(4+\frac{4}{n-2}\right)\left(\frac{(n-1)(n-2)}{2}+\frac{n-2}{2}
\right)>
\frac{5(n-1)(n-2)}{2}
$$
for even $n$, and
$$
\left(4+\frac{4}{n-3}\right)\left(\frac{(n-1)(n-2)}{2}+\frac{n-1}{2}
\right)>
\frac{5(n-1)(n-2)}{2}
$$
for odd $n$. It follows $n\le 9$ which contradicts to the assumption
$n\ge 10$.

This finishes the proof.
\end{proof}

\subsection{Ground fields of arithmetic hyperbolic reflection\\
groups in dimensions $n\ge 6$}\label{subsec:ge6}

Let us introduce one more set of fields. Let us consider plane
(or Fuchsian) arithmetic hyperbolic reflection groups $W$ with a
quadrangle fundamental polygon
$K$ of minimality $14$. We remind that this means that
$P(K)=\{\delta_1,\delta_2,\delta_3,\delta_4\}$ satisfies the condition
$$
\delta_i\cdot \delta_j<14,\ \forall\ \delta_i,\delta_j\in P(K).
$$
Respectively, we call $K$ as {\it arithmetic quadrangle of minimality $14$.}

\begin{definition} We denote by $\Gamma_{2,4}(14)$ the set of
Gram graphs $\Gamma (P(K))$ of all arithmetic quadrangles $K$
of minimality $14$.
The set $\F_{2,4}(14)$ consists of all their ground fields.
\label{F{2,4}(14)}
\end{definition}

By Borel \cite{Borel} and Takeuchi \cite{Tak4}, for fixed $g\ge 0$
and $t\ge 0$,
the number of arithmetic Fuchsian groups with signatures
$(g;\ e_1,\,e_2,\,\dots ,\, e_t)$
is finite. Applying this result to $g=0$ and $t=4$, we obtain that sets of
arithmetic quadrangles $\Gamma_{2,4}$ and their ground fields $\F_{2,4}$
are finite. Then their subsets $\Gamma_{2,4}(14)$ and
$\F_{2,4}(14)$ are also finite.

Moreover, in \cite[pages 383--384]{Tak4} an upper bound $n_0$ of the degree
of ground fields of Fuchsian groups with signatures
$(g;\ e_1,\,e_2,\,\dots,\, e_t)$ is given. It is
\begin{equation}
n_0=(b+\log_e{C(g,t)})/\log_e(a/(2\pi)^{4/3})
\label{Tak1}
\end{equation}
where
$$
a=29.099,\ b=8.3185,\ C(g,t)=2^{2g+t-2}(2g+t-2)^{2/3}
$$
(here $a$ and $b$ are due to Odlyzko).
It follows that
\begin{equation}
[\bk:\bq]\le 11\ \ {\rm for\ \ } \bk\in \F_{2,4}\supset \F_{2,4}(14)
\label{degF2,4(14)}.
\end{equation}

We have the following main result of the paper.

\begin{theorem} In dimensions $n\ge 6$, the ground field of
any arithmetic hyperbolic reflection group belongs to one of finite
sets of fields $\F L^4$, $\F T$, $\F\Gamma^{(4)}_i(14)$, $1\le i\le 5$,
and $\F\Gamma_{2,4}(14)$. In particular, its degree is bounded by $120$.
\label{thfor6}
\end{theorem}

\begin{proof} We use notations of Sect. \ref{subsec:ge10} above.
By \eqref{alpha02}, for $n\ge 6$ a narrow face $\M_e$ has
$\alpha_{n-1}^{(0,2)}<5$.
Thus, $\M_e$ has a triangle or quadrangle (2-dimensional) face.
Let us consider both cases.

By Lemma \ref{lembroken}, we have

\begin{lemma}
Let $W$ be an arithmetic hyperbolic reflection group, $\M$ its fundamental
chamber and their ground field is different from fields from $\F L^4$
and $\F T$.

Then $\M$ has no triangle faces (2-dimensional).
\label{lemtriangle}
\end{lemma}

\begin{proof} Assume $\M$ contains a triangle 2-dimensional face.
Then the edge polyhedron of $\M$ corresponding to the
edge $v_1v_2$ of two vertices $v_1$ and $v_2$ of this
triangle has the corresponding elements $e_1$ and $e_2$
(from Lemma \ref{lembroken})  such that the hyperplanes $\Hh_{e_1}$ and
$\Hh_{e_2}$ have a common point which is the third vertex of the triangle.
Then $e_1\cdot e_2\le 2$.
This contradicts Lemma \ref{lembroken}.

This finishes the proof.
\end{proof}

Let $\M_e$ has a triangle face. By Lemma \ref{lemtriangle}, then the
ground field of $\M$
belongs to $\F L^4$ or $\F T$ as required.

\begin{lemma}
Let $W$ be an arithmetic hyperbolic reflection group, $\M$ its fundamental
chamber and their ground field is different from fields from $\F L^4$, $\F T$
and $\F \Gamma^{(4)}_i$, $1\le i\le 5$.

Let $\M_e$, $e\in P(\M)$,  be a narrow face of $\M$ of minimality $14$.
Let $\M_4$ be a quadrangle face of  $\M_e$.  Let $Q\subset P(\M) $
consists of all $n-2$ elements which
are perpendicular to the plane of $\M_4$, and $\delta_j\in P(\M)$,
$j=1,2,3,4$,  are additional $4$ elements which
are perpendicular to four edges of $\M_4$.

Then all elements $\delta_j$ are perpendicular to $Q$ and $\M_4$ is an
arithmetic quadrangle of minimality $14$
with $P(\M)=\{\delta_1,\delta_2, \delta_3,\delta_4\}$ and with the same
ground field as $\M$. Thus, the ground field of $W$
belongs to $\F_{2,4}(14)$.
\label{lemquadrangle}
\end{lemma}

\begin{proof} We assume that $\delta_1$, $\delta_2$, $\delta_3$ and $\delta_4$
are perpendicular to four consecutive edges of $\M_4$.  The quadrangle
$\M_4$ has a non-right angle.
We can assume that $\delta_1$ and $\delta_2$ are perpendicular to edges
of this angle.  Since the angle is non-right,
$\delta_1$ and $\delta_2$ belong to one connected component of the Coxeter
graph of the vertex of the angle.
By Lemma \ref{lemcoxeter}, then $\delta_1$ and $\delta_2$ give a connected
component of the graph. It follows that $\delta_1$ and
$\delta_2$ are perpendicular to $Q$.

Assume that $\delta_3$ is not perpendicular to $Q$ and $\delta_3\cdot e>0$
for $e\in Q$. Then
$\delta_1,\delta_2,\delta_3, e$ define an edge graph $\Gamma^{(4)}_4$
(if $\delta_2\cdot \delta_3>0$) or
$\Gamma^{(4)}_5$ (if $\delta_2\cdot\delta_3=0$) of minimality $14$.
It follows that the ground field of $\M$ belongs to
$\F \Gamma^{(4)}_4(14)$ or $\F \Gamma^{(4)}_5(14)$, and we get
a contradiction. Thus, $\delta_3$
is perpendicular to $Q$. Similarly we can prove that $\delta_4$ is
perpendicular to $Q$.

This finishes the proof of the lemma.
\end{proof}

If $\M_e$ contains a quadrangle, by  Lemma \ref{lemquadrangle}, the
ground field of $W$ belongs to $\F_{2,4}(14)$ as
required.
This finishes the proof of the theorem. \end{proof}

\subsection{Ground fields of arithmetic hyperbolic reflection\\
groups in dimensions  $n\ge 4$}\label{subsec:ge4}

Unfortunately, in dimension $n\ge 4$ we don't know similar results
to Theorems \ref{thfor10} and \ref{thfor6}.
Possibly, the recent preprint by Agol, Belolipetsky, Storm and Whyte
\cite{ABSW} contains some similar information.
It gives some effective bounds on degrees and discriminants of ground
fields of arithmetic hyperbolic reflection groups for $n\ge 4$.
Unfortunately, they are not explicit, it seems.

\medskip

On the other hand, the following result had been obtained in \cite{Nik5}.

\begin{theorem} (\cite{Nik5}) For $n\ge 4$, the ground field of any
n-dimensional arithmetic hyperbolic reflection group is either
the ground field of one of $n-1$ or $n-2$-dimensional arithmetic hyperbolic
reflection group with a fundamental chamber of minimality $14$,
or a field from
$\F L^4$, $\F T$ and $\F\Gamma^{(4)}_i(14)$, $i=1,\,2,\,3,\,4$.

In particular, its degree is bounded by the maximum of degrees of ground
fields of $n-1$ and $n-2$-dimensional arithmetic hyperbolic reflection
groups with a fundamental chamber of minimality $14$, and of $120$
(according to Theorem \ref{thdegrees} of this paper).
\label{thfor4}
\end{theorem}

\begin{proof} We repeat arguments of \cite{Nik5}.

Let $\M$ be a fundamental chamber of $W$. We can assume that
$\M$ is compact and the ground field of $\M$ is not contained in
$\F L^4$ and $\F T$. Let $\M_e$, $e\in P(\M)$, be a face of $\M$ of
minimality $14$.

If all hyperplanes $\Hh_\delta$, $\delta \in P(\M,e)-\{e\}$, are perpendicular
to $\M_e$ (equivalently, $\delta \cdot e=0$), then
$\M_e$ is a fundamental chamber of
$n-1$-dimensional arithmetic hyperbolic reflection group with the same
ground field as $\M$. Obviously, $P(\M_e)=P(\M,e)-\{e\}$, and $\M_e$ has
minimality $14$.

If this is not the case, there exists $f\in P(\M,e)-\{e\}$
such that $f\cdot e>0$ (equivalently, $f$ and $e$ are connected by a
thin edge in Gram graph of $\M$). Then
$\M_{e,f}=\M\cap \Hh_e\cap \Hh_f$ is $n-2$-dimensional face
of $\M$. Let $P(\M,e,f)$ be the set of all $\delta\in P(\M)$ such that
the hyperplane $\Hh_\delta$ intersects the codimension-two subspace
$\Hh_e\cap \Hh_f$ (then $\M\cap \Hh_e\cap \Hh_f\cap \Hh_\delta$ is a
codimension-three face of $\M$ if $\delta$ is different from $e$ and $f$).
If $\Hh_\delta \perp \Hh_e\cap \Hh_f$ (equivalently,
$\delta\cdot e=\delta\cdot f=0$) for all $\delta\in P(\M,e,f)-\{e,f\}$, then
$\M_{e,f}$ is a fundamental chamber of an arithmetic
hyperbolic reflection group of dimension $n-2$ with the same ground field as
$\M$. Obviously, $P(\M_{e,f})=P(\M,e,f)-\{e,f\}$, and $\M_{e,f}$ has
minimality $14$.

If this is not the case, there is $g\in P(\M,e,f)-\{e,f\}$ such that
$\Hh_g$ is not perpendicular to $\Hh_e\cap \Hh_f$. This means that either
$g\cdot e>0$ or $g\cdot f>0$. Thus, the Gram graph of $e,f,g$ is a
connected negative definite (i. e. connected Coxeter) graph.

We consider an edge $r$ in the face $\M_{e,f}=\M\cap \Hh_e\cap \Hh_f$
of $\M$ such that $r$ terminates in the hyperplane $\Hh_g$.
Thus one of vertices of $r$ is contained in $\Hh_g$ while the other
is not (equivalently, $r$ is not contained in $\Hh_g$).
The existence of such an edge is obvious.
Let $h\in P(\M)$ defines the hyperplane $\Hh_h$ which contains only the
vertex of $r$ which does not belong to $\Hh_g$. Then $g$ and $h$ are joined
by a broken edge in the edge graph $\Gamma(r)$ of $r$ (here we can
additionally assume that the ground field does not belong to
$\F L^4$ and $\F T$). The four elements
$\{e,f,g,h\}$ define a connected hyperbolic subgraph of $\Gamma (r)$
with four vertices. It is one of graphs $\Gamma^{(4)}_i(14)$, $i=1,2,3,4$.
Then the ground field of $\M$ belongs to one of sets
$\F \Gamma^{(4)}_i(14)$, $1\le i\le 4$.
This finishes the proof.
\end{proof}

Theorem \ref{thfor4} shows that ground fields of arithmetic hyperbolic
reflection groups which are different from fields of $\F L^4$, $\F T$ and
$\F \Gamma^{(4)}_i(14)$, $i=1,\,2,\,3,\,4$, come up from $2$-dimensional and
$3$-dimensional arithmetic hyperbolic reflection groups.

\subsection{Ground fields of arithmetic hyperbolic reflection\\
groups in dimension  $n=3$}\label{subsec:deg=3}
Finiteness of the number of
maximal arithmetic hyperbolic reflection groups of
dimension $n=3$ was proved by Agol in \cite{Agol}.
It follows that the number of ground fields of 3-dimensional arithmetic
hyperbolic reflection groups is also finite. Unfortunately,
an explicit  bound of degrees of these fields is not  known (by the
author's knowledge). We remind that by
 \cite{Nik1}, the set of ground fields of arithmetic hyperbolic
reflection groups of the fixed degree is finite and
can be effectively found. Thus, an explicit bound of the degree
is the crucial problem in finding of ground fields of 3-dimensional
hyperbolic reflection groups.

\subsection{Ground fields of arithmetic hyperbolic reflection\\
groups in dimension  $n=2$}\label{subsec:deg=2}
Finiteness of the number of
maximal arithmetic hyperbolic reflection groups of dimension
$n=2$ was proved by Long, Maclachlan and Reid \cite{LMR}. They proved
finiteness of maximal arithmetic Fuchsian groups of genus $0$.
Their ground fields contain ground fields of all arithmetic hyperbolic
reflection groups.

Let $\Gamma$ be a cocompact maximal arithmetic Fuchsian group of genus $0$.
Using results by M.-F. Vigneras \cite{Vig} and Zograf \cite{Zog},
the bound
\begin{equation}
{\rm Area} \left({\mathbb H}^2 /\Gamma\right)\le \frac{128 \pi}{3}
\label{ar2}
\end{equation}
of the area of the arithmetic quotient was obtained in \cite{LMR}.
Let $\Gamma$ has signature $(0; e_1,\dots, e_t)$ where
$e_i\ge 2$.
Then the area of the quotient is equal to
$$
{\rm Area} \left({\mathbb H}^2 /\Gamma\right)=2\pi
\left(t-2-\sum_{i=1}^t{\frac{1}{e_i}}\right).
$$
Since $e_i\ge 2$, by \eqref{ar2}, we obtain $2\pi (t-2-t/2)\le128\pi/3$,
and $t\le 46$. By the result of Takeuchi
\eqref{Tak1}, we obtain $n_0\le 44$.

Thus,{\it degree of the ground field of any arithmetic
Fuchsian group of genus 0 is less or equal to $44$.}
In particular, {\it degree of the ground field of any 2-dimensional
arithmetic hyperbolic reflection group is less or equal to
$44$.}

\medskip

Summarising above results, we see that an explicit bound of degrees of
ground fields of arithmetic hyperbolic reflection groups remains unknown in
dimensions  $n=3,\,4,\,5$ only. Moreover, the dimension $n=3$ is crucial
for this problem. If one finds this bound for $n=3$, we will
know it for all remaining dimensions $n=4$ and $n=5$.


\section{Ground fields of V-arithmetic connected finite edge graphs with
four vertices  of the minimality $14$.}
\label{secfieldedggraphs4}

Here we shall obtain explicit upper bounds of  degrees of fields from
the finite sets $\F \Gamma^{(4)}_i(14)$, $1\le i\le 5$ (see Definition
\ref{fieldsedge4}), and prove Theorem \ref{thdegrees}.  Moreover, our
considerations will deliver important information about these
sets of fields.

Like for the proof of Theorem \ref{tholdVedge} from \cite{Nik2}, we use
the following general results from \cite{Nik2} (we use corrections from Section \ref{sec:Appendix}).

\begin{theorem} (\cite[Theorem 1.2.1]{Nik2}) Let $\bff$ be a totally
real algebraic number field, and let each
embedding $\sigma:\bff\to \br$ corresponds to an interval
$[a_\sigma,b_\sigma]$ in $\br$ where
$$
\prod_{\sigma }{\frac{b_\sigma-a_\sigma}{4}}<1.
$$
In addition, let the natural number $m$ and the intervals
$[s_1,t_1],\dots, [s_m,t_m]$ in $\br$ be fixed. Then there exists
a constant $N(s_i,t_i)$ such that, if $\alpha$ is a totally real
algebraic integer and if the following inequalities hold for the
embeddings $\tau:\bff(\alpha) \to \br$:
$$
s_i\le \tau(\alpha)\le t_i\ \ for\ \ \tau=\tau_1,\dots ,\tau_m,
$$
$$
a_{\tau | \bff}\le \tau(\alpha)\le b_{\tau | \bff}\ \ for\ \
\tau\not=\tau_1,\dots,\tau_m,
$$
then
$$
[\bff(\alpha):\bff]\le N(s_i,t_i).
$$
\label{th121}
\end{theorem}

\begin{theorem} (\cite[Theorem 1.2.2]{Nik2})
Under the conditions of Theorem \ref{th121}, $N(s_i,t_i)$ can be
taken to be $N(s_i,t_i)=N$,
 where $N$ is the least natural number
solution of the inequality
\begin{equation}
N\ln{(1/R)} - M\ln{(2N+2)}-\ln{B}\ge \ln{S}.
\label{cond for n}
\end{equation}

Here
\begin{equation}
M=[\bff : \bq],\ \ \ B=\sqrt{|{\rm discr\ } \bff|};
\label{MB}
\end{equation}
\begin{equation}
R=\sqrt{\prod_\sigma {\frac{b_\sigma-a_\sigma}{4}}},\ \ \
S=\prod_{i=1}^{m}{\frac{2er_i}{b_{\sigma_i}-a_{\sigma_i}}}
\label{RS}
\end{equation}
where
\begin{equation}
\sigma_i=\tau_i|\bff,\ \ \  r_i=\max\{{|b_i-a_{\sigma_i}|,
|b_{\sigma_i}-a_i|}\}.
\label{ri}
\end{equation}
\label{th122}
\end{theorem}

We note that the proof of Theorems \ref{th121} and \ref{th122}
uses a variant of Fekete's Theorem (1923) about existence of
non-zero integer polynomials of bounded degree which differ only
slightly from zero on appropriate intervals. See \cite[Theorem
1.1.1]{Nik2} (see its corrections in Section \ref{sec:Appendix}, Theorems 
\ref{thFekete1}, \ref{thFekete2}).

\subsection{Fields from $\F \Gamma_1^{(4)}(14)$}\label{subsec:gam(4)1}

\begin{figure}
\begin{center}
\includegraphics[width=6cm]{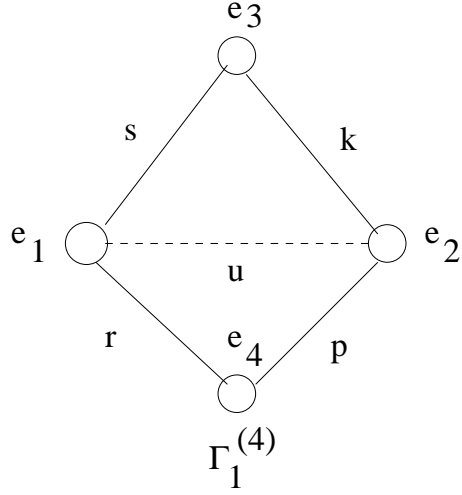}
\end{center}
\caption{The graph $\Gamma_1^{(4)}$}
\label{graphg41}
\end{figure}

For $\Gamma_1^{(4)}(14)$ (see Figure \ref{graphg41}) we assume that integers
$s,k,r,p\ge 3$. Subgraphs
$\Gamma_1^{(4)}-\{e_1\}$ and $\Gamma_1^{(4)}-\{e_2\}$ must be Coxeter
graphs. It follows that
we must consider only (up to obvious symmetries) the following cases:
either $s=k=3$ and $5\ge r\ge p\ge 3$, or $s=p=3$ and $5\ge r\ge k\ge 4$;
the totally real algebraic integer $u$ satisfies the inequality
$2<u<14$.

The corresponding Gram matrix is
\begin{equation}
\left(\begin{array}{cccc}
-2                            & u
& 2\cos{\frac{\pi}{s}}  & 2\cos{\frac{\pi}{r}}\\
u                              & -2
& 2\cos{\frac{\pi}{k}}  & 2\cos{\frac{\pi}{p}}\\
2\cos{\frac{\pi}{s}} & 2\cos{\frac{\pi}{k}} & -2
& 0                             \\
2\cos{\frac{\pi}{r}}  & 2\cos{\frac{\pi}{p}} & 0
& -2
\end{array}\right) \ ,
\label{gramm1(4)}
\end{equation}
where all entries are algebraic integers,
and its determinant $d(u)$ is given by the equality
\begin{equation}
\begin{array}{l}
-\frac{d(u)}{4}=\\
\left(u+2\left(\cos{\frac{\pi}{r}}\cos{\frac{\pi}{p}} +
\cos{\frac{\pi}{k}}\cos{\frac{\pi}{s}}\right)\right)^2
-\left(\cos{\frac{2\pi}{k}}+\cos{\frac{2\pi}{p}}\right)
\left(\cos{\frac{2\pi}{r}}+\cos{\frac{2\pi}{s}}\right) .
\end{array}
\label{detgramm1(4)}
\end{equation}

The ground field $\bk$ is generated by the cyclic products
$$
\bk=(u^2, u\cos{(\pi/s)}\cos{(\pi/k)},u\cos{(\pi/p)}\cos{(\pi/r)},
\cos^2{(\pi/s)}, \cos^2{(\pi/k)},
$$
$$
 \cos^2{(\pi/p)}, \cos^2{(\pi/r)},
\cos{(\pi/s)}\cos{(\pi/k)}\cos{(\pi/p)}\cos{(\pi/r)})\,.
$$

Here and in what follows we always denote by $\sigma^{(+)}:\bk\to
\br$ the geometric (the identity) embedding, and by $\sigma:\bk
\to \br$ all other embeddings $\sigma\not=\sigma^{(+)}$. We have
$\sigma^{(+)}(d(u))<0$ because $\Gamma^{(4)}_1$ is hyperbolic for
$\sigma^{(+)}$ since $e_1, e_2$ define a hyperbolic subgraph
because $u=e_1\cdot e_2>2$. And $\sigma(d(u))>0$ since
$\Gamma^{(4)}_1$ is negative definite for $\sigma$. In particular,
$-2<\sigma (u)<2$ by Cauchy inequality. Thus, $\Gamma^{(4)}_1$ is
V-arithmetic if and only if
\begin{equation}
\sigma\left(u+2(\cos{\frac{\pi}{r}}\cos{\frac{\pi}{p}} +
\cos{\frac{\pi}{k}}\cos{\frac{\pi}{s}})\right)^2<
$$
$$
\sigma\left((\cos{\frac{2\pi}{k}}+\cos{\frac{2\pi}{p}})
(\cos{\frac{2\pi}{r}}+\cos{\frac{2\pi}{s}})\right).
\label{V1(4)}
\end{equation}
for each $\sigma\not=\sigma^{(+)}$.

We have $\bk=\bq(u^2)$ since $4<\sigma^{(+)}(u^2)<14^2$ and
$0<\sigma(u^2)<4$ if $\sigma\not=\sigma^{(+)}$.
We have $[\bk(u):\bk]=2$ if $u\notin \bk$.
If $\tau:\bk(u) \to \br$ gives $\tau|\bk=\sigma^{(+)}$,
then $\tau(u)=\pm u$ (where $u$ is taken for
the geometric embedding $\sigma^{(+)}$), and  either
$2<\tau(u)<14$ or $-14<\tau(u)<-2$.
The last inequality is possible, only if $u$ does not belong to $\bk$.

If $\tau|\bk=\sigma\not=\sigma^{(+)}$, then by \eqref{V1(4)},
$$
-2\widetilde\tau\left(\cos{\frac{\pi}{r}}\cos{\frac{\pi}{p}} +
\cos{\frac{\pi}{k}}\cos{\frac{\pi}{s}}\right)-
$$
$$
-\sqrt{\tau\left((\cos{\frac{2\pi}{k}}+\cos{\frac{2\pi}{p}})
(\cos{\frac{2\pi}{r}}+\cos{\frac{2\pi}{s}})\right)}.
$$
$$
<\tau(u)<-2\widetilde \tau\left(\cos{\frac{\pi}{r}}\cos{\frac{\pi}{p}} +
\cos{\frac{\pi}{k}}\cos{\frac{\pi}{s}}\right)+
$$
$$
\sqrt{\tau\left((\cos{\frac{2\pi}{k}}+\cos{\frac{2\pi}{p}})
(\cos{\frac{2\pi}{r}}+\cos{\frac{2\pi}{s}})\right)}
$$
where $\tilde\tau$ extends $\tau$.
Thus, $\tau(u)$ belongs to the interval of the length
$$
2\sqrt{\tau\left((\cos{\frac{2\pi}{k}}+\cos{\frac{2\pi}{p}})
(\cos{\frac{2\pi}{r}}+\cos{\frac{2\pi}{s}})\right)}.
$$
We apply Theorems \ref{th121} and \ref{th122} to
$$
\bff=\bq( \cos^2{(\pi/s)}, \cos^2{(\pi/k)}, \cos^2{(\pi/p)},
\cos^2{(\pi/r)})
$$
and $\alpha=u$
to bound $[\bff(u):\bff]$. Since $\bk=\bff(u^2)$,
then $m=[\bff(u):\bk]\le 2$.
From considerations above, we have
$$
M=[\bff:\bq],\ \ B=\sqrt{|{\rm discr\ }\bff |},
$$
$$
R=\frac{N_{\bff/\bq}\left((\cos{\frac{2\pi}{k}}+\cos{\frac{2\pi}{p}})
(\cos{\frac{2\pi}{r}}+\cos{\frac{2\pi}{s}})\right)^{1/4}}
{2^{M/2}},
$$
$$
S=\left(\frac{16\cdot e}
{\sqrt{\left((\cos{\frac{2\pi}{k}}+\cos{\frac{2\pi}{p}})
(\cos{\frac{2\pi}{r}}+
\cos{\frac{2\pi}{s}})\right)}}\right)^m,
$$
$$
[\bk:\bq]\le NM/m
$$
where $N$ is the least natural solution of the inequality
$$
N\ln{(1/R)} - M\ln{(2N+2)}-\ln{B}\ge \ln{S}.
$$

By direct calculations, we obtain the following.

Let $k=s=3$. Then $m=1$.

For $r=p=3$, we obtain $M=1$, $B=1$, $R=1/\sqrt{2}$, $S=16e$.
Then $N=22$ and $[\bk:\bq]\le 22$.

For $r=4$, $p=3$, we obtain  $M=1$, $B=1$, $R=1/2^{3/4}$, $S=16e\sqrt{2}$.
Then $N=15$ and $[\bk:\bq]\le 15$.

For $r=5$, $p=3$, we obtain $M=2$, $B=\sqrt{5}$, $R=1/2^{3/2}$, \linebreak
$S=16e/\sqrt{1/2-cos(2\pi/5)}$. Then $N=12$ and $[\bk:\bq]\le 24$.

For $r=p=4$,  we obtain $M=1$, $B=1$, $R=1/2$, $S=32e$. Then $N=12$
and $[\bk:\bq]\le 12$.

For $r=5$, $p=4$, we obtain $M=2$, $B=\sqrt{5}$, $R=1/4$,\linebreak
$S=32e/\sqrt{1-2cos(2\pi/5)}$. Then $N=9$ and $[\bk:\bq]\le 18$.

For $r=p=5$, we obtain $M=2$, $B=\sqrt{5}$, $R=1/4$,
$S=16e/(1/2-cos(2\pi/5))$. Then $N=9$ and $[\bk:\bq]\le 18$.

Now let $s=p=3$. Then $m=1$ or $m=2$.

For $r=k=4$, we obtain $M=1$, $B=1$, $R=1/2$, $S=(32 e)^m$.
Then $N=12$ and $[\bk:\bq]\le 12$ for $m=1$; $N=19$ and
$[\bk:\bq]\le 9$ for $m=2$.

For $r=5$, $k=4$, we obtain $m=1$, $M=2$, $B=\sqrt{5}$, $R=1/4$,
$S=32e/\sqrt{1-2cos(2\pi/5)}$. Then $N=9$ and $[\bk:\bq]\le 18$.

For $r=k=5$, we obtain $M=2$, $B=\sqrt{5}$, $R=1/4$,
$S=(16e/(1/2-cos(2\pi/5))^m$. Then $N=9$ and $[\bk:\bq]\le 18$ for
$m=1$, and $N=14$ and $[\bk:\bq]\le 14$ for $m=2$.

Note that the bound for $[\bk:\bq]$ is always worse for $m=1$ than for $m=2$.
It follows from our method. Further, in  similar considerations, we
can consider $m=1$ only.

Thus, our upper bound for degrees of fields from
$\F \Gamma^{(4)}_1(14)$ is $24$.

\subsection{Fields from $\F \Gamma_2^{(4)}(14)$}\label{subsec:gam(4)2}
For $\Gamma_2^{(4)}(14)$ (see Figure \ref{graphg42}),
$s,k,p \ge 3$ are natural numbers and $2<u<14$ is
a totally real algebraic integer. Moreover, we have only the following
possibilities: $3\le s\le k\le 5$, $p=3$; $s=k=3$, $p=4,\,5$.

The ground field $\bk=\bq(u^2)$ contains cyclic products
$$
\cos^2{\frac{\pi}{s}},\ \cos^2{\frac{\pi}{k}},\ \cos^2{\frac{\pi}{p}},\
u^2,\ u\,\cos{\frac{\pi}{s}}\cos{\frac{\pi}{k}}\,.
$$
This case is similar to $\F\Gamma^{(4)}_1(14)$. The determinant $d(u)$
of the Gram matrix is determined by the equality
$$
-\frac{d(u)}{4}=
\sin^2{\frac{\pi}{p}}\,u^2+4\cos{\frac{\pi}{s}} \cos{\frac{\pi}{k}}\,u+
4 (\cos^2{\frac{\pi}{s}}+\cos^2{\frac{\pi}{k}}+\cos^2{\frac{\pi}{p}}-1).
$$
Let
$$
D=16\cos^2{\frac{\pi}{s}} \cos^2{\frac{\pi}{k}}+
16 \sin^2{\frac{\pi}{p}}(1-\cos^2{\frac{\pi}{s}}-\cos^2{\frac{\pi}{k}}-
 \cos^2{\frac{\pi}{p}})
$$
be the discriminant of this quadratic polynomial of the variable $u$.
The graph $\Gamma^{(4)}_4(14)$ is V-arithmetic if and only if for
$\tau:\bk(u)\to \br$ which is different from $\sigma^{(+)}$ on $\bk$,
one has
$$
\frac{-4\widetilde\tau\left(\cos{\frac{\pi}{s}} \cos{\frac{\pi}{k}}\right)-
\sqrt{\tau(D)}}{2\tau(\sin^2{\frac{\pi}{p}})}\
<\tau(u)<\
$$
$$
< \frac{-4\widetilde\tau\left(\cos{\frac{\pi}{s}} \cos{\frac{\pi}{k}}\right)+
\sqrt{\tau(D)}}{2\tau(\sin^2{\frac{\pi}{p}})}
$$
where $\widetilde{\tau}$ extends $\tau$.

\begin{figure}
\begin{center}
\includegraphics[width=6cm]{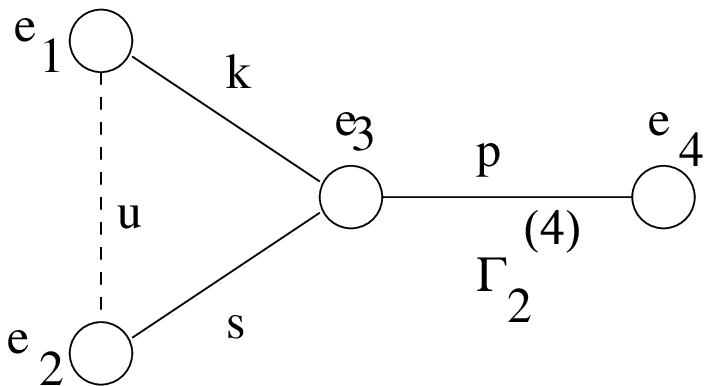}
\end{center}
\caption{The graph $\Gamma_2^{(4)}$}
\label{graphg42}
\end{figure}

Thus, $\tau(u)$ belongs to an interval of the length
$2\sqrt{\tau\left(D/(4\sin^4{\frac{\pi}{p}})\right)}$.

We can apply Theorems \ref{th121}, \ref{th122} to
$\bff=\bq(\cos^2{\frac{\pi}{s}},\ \cos^2{\frac{\pi}{k}},\
\cos^2{\frac{\pi}{p}})$ and
$\alpha=u$. Then
$$
M=[\bff:\bq],\ \ B=\sqrt{|{\rm discr\ }\bff |},
$$
$$
R=\frac{{N_{\bff/\bq}(D)}^{1/4}}
{N_{\bff/\bq}\left( \sin^2{(\pi/p)}\right)^{1/2}2^M}\ ,\ \
S=\frac{32e\sin^2{(\pi/p)}}
{\sqrt{D}},
$$
$[\bk:\bq]\le NM$ where $N$ is the smallest natural solution of the inequality
$$
N\ln{(1/R)} - \ln{(2N+2)}-\ln{B}\ge \ln{S}.
$$

We obtain:

\noindent
if $s=k=p=3$,  then $[\bk:\bq]\le 39$;
if $s=3$, $k=4$, $p=3$, then $[\bk:\bq]\le 21$;

\noindent
if $s=3$, $k=5$, $p=3$, then $[\bk:\bq]\le 34$;
if $s=k=4$, $p=3$, then $[\bk:\bq]\le 14$;

\noindent
if $s=4$, $k=5$, $p=3$, then $[\bk:\bq]\le 22$;
if $s=k=5$, $p=3$, then $[\bk:\bq]\le 24$;

\noindent
if $s=k=3$, $p=4$, then $[\bk:\bq]\le 22$;
if $s=k=3$, $p=5$, then $[\bk:\bq]\le 32$.

Thus, our upper bound for degrees of fields from
$\F \Gamma^{(4)}_2(14)$ is $39$.

\subsection{Fields from $\F \Gamma_3^{(4)}(14)$}\label{subsec:gam(3)4}
For $\Gamma_3^{(4)}(14)$ (see Figure \ref{graphg43}),
$s\ge 2$, $k,r \ge 3$ are natural numbers, and $2<u<14$ is
a totally real algebraic integer. Moreover, we have only the following
possibilities: $s=2$, $k=3$, $r=3,\,4,\,5$;
$s=2$, $k=4,\,5$, $r=3$; $3\le s\le k\le 5$, $r=3$; $s=k=3$, $r=4,\,5$.

\begin{figure}
\begin{center}
\includegraphics[width=6cm]{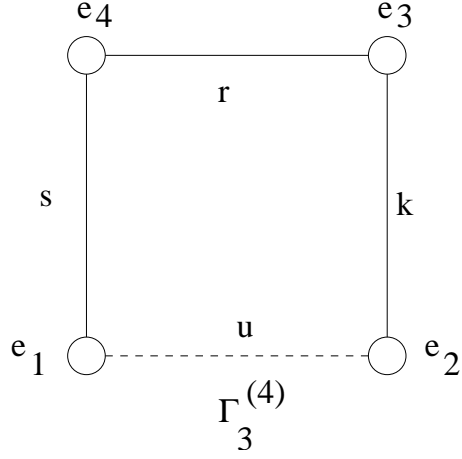}
\end{center}
\caption{The graph $\Gamma_3^{(4)}$}
\label{graphg43}
\end{figure}

The ground field $\bk=\bq(u^2)$ contains cyclic products
$$
\cos^2{\frac{\pi}{s}},\ \cos^2{\frac{\pi}{k}},\ \cos^2{\frac{\pi}{r}},\
u^2,\ u\,\cos{\frac{\pi}{s}}\cos{\frac{\pi}{k}}\cos{\frac{\pi}{r}}\,.
$$

This case is similar to $\F\Gamma^{(4)}_1(14)$ and $\F\Gamma^{(4)}_2(14)$.
The determinant $d(u)$  of the Gram matrix is determined by the equality
$$
-\frac{d(u)}{4}=
\sin^2{\frac{\pi}{r}}\,u^2+2\cos{\frac{\pi}{s}} \cos{\frac{\pi}{k}}
\cos{\frac{\pi}{r}}\,u+4 \cos^2{\frac{\pi}{r}}-
4 \sin^2{\frac{\pi}{s}}\sin^2{\frac{\pi}{k}}.
$$
Let
$$
D=4\cos^2{\frac{\pi}{s}} \cos^2{\frac{\pi}{k}} \cos^2{\frac{\pi}{r}}+
16 \sin^2{\frac{\pi}{s}}\sin^2{\frac{\pi}{k}}\sin^2{\frac{\pi}{r}}-
16 \sin^2{\frac{\pi}{r}}\cos^2{\frac{\pi}{r}}
$$
be the discriminant of this quadratic polynomial of the variable $u$.
The graph $\Gamma^{(4)}_3(14)$ is V-arithmetic if and only if for
$\tau:\bk(u)\to \br$ which is different from $\sigma^{(+)}$ on $\bk$,
one has
$$
\frac{-2\widetilde\tau\left(\cos{\frac{\pi}{s}} \cos{\frac{\pi}{k}}
\cos{\frac{\pi}{r}}\right)-\sqrt{\tau(D)}}{2\tau(\sin^2{\frac{\pi}{r}})}\
<\tau(u)<\
$$
$$
< \frac{-2\widetilde\tau\left(\cos{\frac{\pi}{s}} \cos{\frac{\pi}{k}}
\cos{\frac{\pi}{r}}\right)+\sqrt{\tau(D)}}{2\tau(\sin^2{\frac{\pi}{r}})}
$$
where $\widetilde{\tau}$ extends $\tau$.

Thus, $\tau(u)$ belongs to an interval of the length
$2\sqrt{\tau\left(D/(4\sin^4{\frac{\pi}{r}})\right)}$

We can apply Theorems \ref{th121}, \ref{th122} to
$\bff=\bq(\cos^2{\frac{\pi}{s}},\ \cos^2{\frac{\pi}{k}},\
\cos^2{\frac{\pi}{r}})$ and
$\alpha=u$. Then
$$
M=[\bff:\bq],\ \ B=\sqrt{|{\rm discr\ }\bff |},
$$
$$
R=\frac{{N_{\bff/\bq}(D)}^{1/4}}
{N_{\bff/\bq}\left( \sin^2{(\pi/r)}\right)^{1/2}2^M}\ ,\ \
S=\frac{32e\sin^2{(\pi/r)}}
{\sqrt{D}},
$$
$[\bk:\bq]\le NM$ where $N$ is the smallest natural solution of
the inequality
$$
N\ln{(1/R)} - M\ln{(2N+2)}-\ln{B}\ge \ln{S}.
$$

We obtain:

\noindent
if $s=2$, $k=r=3$, then $[\bk:\bq]\le 83$;
if $s=2$, $k=3$, $r=4$, then $[\bk:\bq]\le 45$;

\noindent
if $s=2$, $k=3$, $r=5$, then $[\bk:\bq]\le 66$;
if $s=2$, $k=4$, $r=3$, then $[\bk:\bq]\le 28$;

\noindent
if $s=2$, $k=5$, $r=3$, then $[\bk:\bq]\le 48$;
if $s=k=r=3$, then $[\bk:\bq]\le 37$;

\noindent
if $s=r=3$, $k=4$, then $[\bk:\bq]\le 18$;
if $s=r=3$, $k=5$, then $[\bk:\bq]\le 24$;

\noindent
if $s=4$, $k=4$, $r=3$, then $[\bk:\bq]\le 9$;
if $s=4$, $k=5$, $r=3$, then $[\bk:\bq]=2$;

\noindent
if $s=k=5$, $r=3$, then $[\bk:\bq]=2$;
if $s=k=3$, $r=4$, then $[\bk:\bq]\le 17$;

\noindent
if $s=k=3$, $r=5$, then $[\bk:\bq]=2$.

For $s=2$ we can improve these estimates considering $\alpha=u^2$.
In this case,
$0<\tau(u^2)<\tau\left(D/(4\sin^4{(\pi/r)})\right)$, and we can
apply Theorems \ref{th121}, \ref{th122} to the same $\bff$
 and
$$
R=\frac{N_{\bff/\bq}(D)^{1/2}}{N_{\bff/\bq}(\sin^2{(\pi/r)})4^M},\ \
S=\frac{2\cdot e\cdot 14^2\cdot 4\cdot \sin^4{(\pi/r)}}{D}\ .
$$
We obtain:

\noindent
if $s=2$, $k=r=3$, then $[\bk:\bq]\le 53$;
if $s=2$, $k=3$, $r=4$, then $[\bk:\bq]\le 31$;

\noindent
if $s=2$, $k=3$, $r=5$, then $[\bk:\bq]\le 32$;
if $s=2$, $k=4$, $r=3$, then $[\bk:\bq]\le 19$;
if $s=2$, $k=5$, $r=3$, then $[\bk:\bq]\le 22$.

Thus, our upper bound for degrees of fields from
$\F \Gamma^{(4)}_3(14)$ is $53$.

\subsection{Fields from  $\F \Gamma_4^{(4)}(14)$}\label{subsec:gam(4)4}
For $\Gamma_4^{(4)}(14)$ (see Figure \ref{graphg44}),
$k\ge 2$, $s,\,r \ge 3$ are natural numbers and $2<u<14$ is
a totally real algebraic integer. Moreover, we have only the following
possibilities: $s=3$, $r=3,\,4,\,5$; $s=4,5$, $r=3$.

\begin{figure}
\begin{center}
\includegraphics[width=6cm]{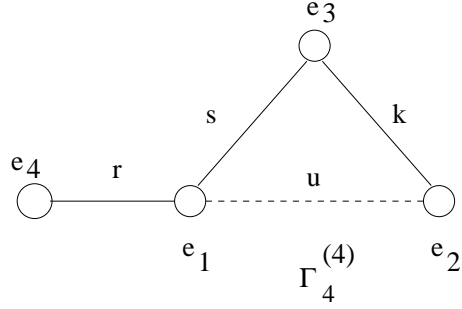}
\end{center}
\caption{The graph $\Gamma_4^{(4)}$}
\label{graphg44}
\end{figure}

The ground field $\bk=\bq(u^2)$ contains cyclic products
$$
\cos^2{\frac{\pi}{s}},\ \cos^2{\frac{\pi}{r}},\ \cos^2{\frac{\pi}{k}},\
u^2,\ u\,\cos{\frac{\pi}{s}}\cos{\frac{\pi}{k}}\,.
$$
This case can be considered as a specialization of the graph
$\Gamma^{(4)}_1$ when we take $p=2$. The determinant $d(u)$ of the
Gram matrix is determined by the equality
$$
\frac{-d(u)}{4}=u^2+4\cos{\frac{\pi}{s}}\cos{\frac{\pi}{k}}\,u+4\,
\cos^2{\frac{\pi}{s}}-4\sin^2 {\frac{\pi}{k}}\sin^2{\frac{\pi}{r}}.
$$
Here
$$
\widetilde{u}=u^2+4\cos{\frac{\pi}{s}}\cos{\frac{\pi}{k}}\,
u+4\cos^2{\frac{\pi}{s}}
$$
is a totally positive algebraic integer
(since minimum of this quadratic polynomial of $u$ is equal to
$4\cos^2{\frac{\pi}{s}}\sin^2{\frac{\pi}{k}}$) which belongs to $\bk$,
and $\Gamma^{(4)}_4$
is V-arithmetic if and only if
$$
0<\sigma(4\cos^2{\frac{\pi}{s}}\sin^2{\frac{\pi}{k}})\le
\sigma(\widetilde{u})<\sigma(4\,\sin^2{\frac{\pi}{r}}\sin^2{\frac{\pi}{k}})
<4
$$
for any $\sigma:\bk\to \br$ which is different from identity $\sigma^{(+)}$.
For $\sigma^{(+)}$, we have
$$
4<\sigma^{(+)}(\widetilde{u})< 14^2+4\cdot 14+4=16^2.
$$
It follows that $\bk=\bq(\widetilde{u})$. We can apply Theorems
\ref{th121} and \ref{th122} to $\bff=\bq(\cos^2{\frac{\pi}{s}},$
$\cos^2{\frac{\pi}{r}},\ \cos^2{\frac{\pi}{k}})$ and
$\alpha=\widetilde{u}$ to estimate $[\bk:\bq]$. We can take
$M=[\bff:\bq]$, $B=\sqrt{{\rm |discr\ }\bff| }$,
$$
R=\sqrt{N_{\bff/\bq}\left((\sin^2{\frac{\pi}{r}}-\cos^2{\frac{\pi}{s}})
\sin^2{\frac{\pi}{k}}\right)},\ \
S=\frac{2\cdot 16^2\cdot e}
{4(\sin^2{\frac{\pi}{r}}-\cos^2{\frac{\pi}{s}})
\sin^2{\frac{\pi}{k}}}\,.
$$
Then $[\bk:\bq]\le MN$ where $N$ is the least natural solution
of the inequality
$$
N\ln{(1/R)} - M\ln{(2N+2)}-\ln{B}\ge \ln{S}.
$$
It follows that $[\bk:\bq]\le 31$ for $2\le k\le 6$ with the maximum
bound $31$ for $k=2$ and $s=r=3$.

When $k\ge 7$, additionally, we should use the
inequality
\begin{equation}
\frac{16^2N_{\bk/\bq}(4\,\sin^2{\frac{\pi}{r}}\sin^2{\frac{\pi}{k}})}
{4\,\sin^2{\frac{\pi}{r}}\sin^2{\frac{\pi}{k}}}>
N_{\bk/\bq}(\widetilde{u})\ge 1
\label{eq42N}
\end{equation}
which follows from considerations above. This additional arguments
will be very similar to much more difficult case of $\Gamma_5^{(4)}(14)$
which we will consider below.

Our upper bound for degrees of fields from
$\F \Gamma^{(4)}_4(14)$ will be $120$ (look at the end of the next section
\ref{subsec:gam(4)5}).

\subsection{Fields from $\F \Gamma_5^{(4)}(14)$}\label{subsec:gam(4)5}
For $\Gamma_5^{(4)}(14)$ (see Figure \ref{graphg45}),
$k\ge s\ge 3$ are natural numbers and $2<u<14$ is a totally real
algebraic integer.

The ground field $\bk=\bq(u^2)$ contains cyclic products
$$
\cos^2{\frac{\pi}{k}},\ \cos^2{\frac{\pi}{s}},\ u^2\ .
$$
The determinant $d(u)$ of the Gram matrix is determined
by the equality
$$
-\frac{d(u)}{4}=u^2-4\sin^2 {\frac{\pi}{k}}\sin^2{\frac{\pi}{s}}.
$$
The $\Gamma^{(4)}_5$ is V-arithmetic if and only if
$$
\sigma(u^2)<\sigma(4\,\sin^2{\frac{\pi}{k}}\sin^2{\frac{\pi}{s}})<4
$$
for any $\sigma:\bk\to \br$ which is different from identity $\sigma^{(+)}$.
For $\sigma^{(+)}$, we have
$$
4<\sigma^{(+)}(u^2)< 14^2.
$$

\begin{figure}
\begin{center}
\includegraphics[width=6cm]{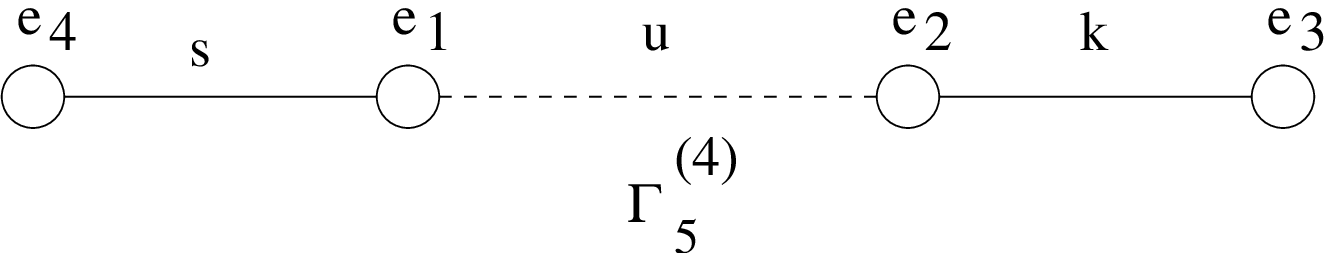}
\end{center}
\caption{The graph $\Gamma_5^{(4)}$}
\label{graphg45}
\end{figure}

We can apply Theorems
\ref{th121} and \ref{th122} to $\bff=\bq(\cos^2{\frac{\pi}{k}},\
\cos^2{\frac{\pi}{s}})$ and
$\alpha=u^2$ to estimate $[\bk:\bq]$. We can take
$M=[\bff:\bq]$, $B=\sqrt{|{\rm discr\ }\bff |}$,
$$
R=\sqrt{N_{\bff/\bq}\left(\sin^2{\frac{\pi}{k}}
\sin^2{\frac{\pi}{s}}\right)}\ ,
\ \
S=\frac{14^2\cdot e}
{2(\sin^2{\frac{\pi}{k}}\sin^2{\frac{\pi}{s}})}.
$$
Then $[\bk:\bq]\le MN_1$ where $N_1$ is the least natural solution
of the inequality
$$
N_1\ln{(1/R)} - M\ln{(2N_1+2)}-\ln{B}\ge \ln{S}.
$$
For each fixed pair  $k\ge s\ge 3$
we can do it and obtain an estimate of the degree $[\bk:\bq]$.
Let us call this method as the {\it Method A.}

From our considerations above, we obtain the inequality
\begin{equation}
\frac{14^2N_{\bk/\bq}(4\,\sin^2{\frac{\pi}{k}}\sin^2{\frac{\pi}{s}})}
{4\,\sin^2{\frac{\pi}{k}}\sin^2{\frac{\pi}{s}}}>
N_{\bk/\bq}(u^2)\ge 1\,.
\label{eq52N}
\end{equation}

We use the following elementary facts about cyclotomic fields. Let
$\bff_l=\bq(\cos^2{(\pi/l)})$ where $l\ge 3$. Then
\begin{equation}
[\bff_l:\bq]=\varphi(l)/2
\label{eqfl:q}
\end{equation}
where $\varphi$ is the Euler function. Let $\bff_{l,m}=\bq(\cos^2{(\pi/l)},
\cos^2{(\pi/m)})$ where $l,m \ge 3$. Then
\begin{equation}
[\bff_{l,m}:\bq]=\frac{\varphi([l,m])}{2\rho(l,m)}
\label{degflm}
\end{equation}
where $[\ ,\ ]$ denotes the least common multiple and $\rho(l,m)=2$ if
$(l,m)|2$, and $\rho(l,m)=1$ otherwise. Here $(\ ,\ )$ denotes
the greatest common divisor.

We have for $l\ge 3$
\begin{equation}
N_{\bff_l/\bq }(4\sin^2{(\pi/l)})=\gamma(l)=\left\{
\begin{array}{cl}
p &\ {\rm if}\  l=p^t>2\   {\rm where\ } p\  {\rm is\  prime,} \\
1 &\ {\rm otherwise.}
\end{array}\
\right .
\label{defgam}
\end{equation}

We denote $\bff=\bff_{k,s}$, $N=[\bk:\bq]$, $n=[k,s]$,
$[\bk:\bff_{k,s}]=m$. By  \eqref{defgam}, we have
$$
N_{\bk/\bq}(4\,\sin^2{\frac{\pi}{k}}\sin^2{\frac{\pi}{s}})=
\frac{N_{\bk/\bq}(4\,\sin^2{\frac{\pi}{k}})N_{\bk/\bq}(4\,
\sin^2{\frac{\pi}{s}})}
{N_{\bk/\bq}(4)}=
$$
$$
\frac{\gamma(k)^{2N/\varphi(k)}\gamma(s)^{2N/\varphi(s)}}{4^N}=
\left(\frac{\gamma(k)^{2/\varphi(k)}\gamma(s)^{2/\varphi(s)}}{4}\right)^N
$$
where
$(\varphi(n)/2\rho(k,s))|N$.

Hence, by \eqref{eq52N},  \eqref{eqfl:q}, \eqref{degflm}, \eqref{defgam},
we obtain
\begin{equation}
N\left(
\ln{2}-\frac{\ln{\gamma(k)}}{\varphi(k)}-\frac{\ln{\gamma(s)}}
{\varphi(s)}
\right)<\\
\ln{7}-\ln{\sin{\frac{\pi}{k}}}-\ln{\sin{\frac{\pi}{s}}},\ \ \
\varphi(n)/2\rho(k,s)| N.
\label{maininN}
\end{equation}

By exact formulae for $\gamma(k)$ and $\varphi(k)$,
it is easy to prove that there exists only finite number of
{\it exceptional pairs} $(k,s)$
such that $k\ge s\ge 3$ and
\begin{equation}
\ln{2}-\frac{\ln{\gamma(k)}}{\varphi(k)}-\frac{\ln{\gamma(s)}}
{\varphi(s)}\le 0.
\label{excpairs}
\end{equation}
They are $s=3$, $k=3, 4, 5, 7, 8, 9, 11,13,17, 19$;
$s=4$, $k=4,5$; $s=5$, $k=5, 7$.

For non-exceptional $k\ge s$ we get
\begin{equation}
\frac{\ln{7}-\ln{\sin{\frac{\pi}{k}}}-\ln{\sin{\frac{\pi}{s}}}}
{ \ln{2}-\frac{\ln{\gamma(k)}}{\varphi(k)}-\frac{\ln{\gamma(s)}}
{\varphi(s)}
}>\frac{\varphi([k,s])}{2\rho(k,s)}, \ \
\label{boundforks}
\end{equation}

\begin{equation}
[\bk:\bff_{k,s}]\le \left[\frac{\ln{7}-\ln{\sin{\frac{\pi}{k}}}-
\ln{\sin{\frac{\pi}{s}}}}
{ \ln{2}-\frac{\ln{\gamma(k)}}{\varphi(k)}-\frac{\ln{\gamma(s)}}
{\varphi(s)}
}\left/\frac{\varphi([k,s])}{2\rho(k,s)} \right.\right],\,
\label{boundforKF}
\end{equation}
and
\begin{equation}
N=[\bk:\bq]\le \left[\frac{\ln{7}-\ln{\sin{\frac{\pi}{k}}}-
\ln{\sin{\frac{\pi}{s}}}}
{ \ln{2}-\frac{\ln{\gamma(k)}}{\varphi(k)}-\frac{\ln{\gamma(s)}}
{\varphi(s)}
}\left/\frac{\varphi([k,s])}{2\rho(k,s)} \right.\right]\cdot
\frac{\varphi([k,s])}{2\rho(k,s)}\ .
\label{boundforK}
\end{equation}

Let us show that there exists only a finite number of non-exceptional
pairs $k\ge s\ge 3$ satisfying the inequality \eqref{boundforks}.

By exact formulae for $\gamma(l)$ and $\varphi(l)$, it is easy to find the
minimum of
$\ln{2}-\ln{\gamma(k)}/\varphi(k)-\ln{\gamma(s)}/\varphi(s)$
for all non-exceptional pairs. It is achieved for $k=23$, $s=3$, and it
is equal to $ln{2}-log(23)/22-log(3)/2\ge 0.00131857$.
Since $x\ge \sin{x}$ for small $x$, $\rho(k,s)\le 2$ and
$\varphi([k,s])\ge \varphi(k)$, we obtain
$ln{7}-2ln{\pi}+2ln(k)\ge 0.00131857\cdot \varphi(k)/4$. Using the trivial
estimate $\varphi(k)\ge \sqrt{k-2}$, we obtain
$ln{7}-2ln{\pi}+2ln(k)\ge 0.0003296425\sqrt{k-2}$.
It follows that $s\le k < 2.1\cdot 10^{10}$. It follows the finiteness.

Using a better estimate $\varphi(k)\ge Ck/ln(ln(k))$,
one can get a better estimate for $k$. One can take
$C=\varphi(6)log(log(6))/6\ge 0.19439$ for $k\ge 6$. See \cite{Vinog}.
It follows, $s\le k < 10^7$.

It follows, that all non-exceptional pairs $k\ge s$
satisfying \eqref{boundforks} can be found
(using a computer). Using \eqref{boundforKF} and \eqref{boundforK},
the bounds for $[\bk:\bff_{k,s}]$ and for $N=[\bk:\bq]=
[\bk:\bff_{k,s}][\bff_{k,s}:\bq]$ can be found for each such a pair.
This we call the {\it Method B}.

If \eqref{boundforK} gives a poor bound for $N$ because of the
bound \eqref{boundforKF} for $[\bk:\bff_{k,s}]$ is poor,
we can improve the bound for $[\bk:\bff_{k,s}]$ using the Method A above.
Also we can apply the Method A to all exceptional pairs $k\ge s$.

As a result, we obtain that $[\bk:\bq]\le 120$. This is achieved
for $k=31$, $s=3$. In this case, \eqref{boundforK} gives
$[\bk:\bq]\le 165$, and \eqref{boundforKF} gives  $[\bk:\bff_{31,3}]\le 11$.
But, the Method A improves the last estimate and gives $[\bk:\bff_{31,3}]\le 8$. 
Thus, we obtain $[\bk:\bq]\le 15\cdot 8=120$ since $[\bff_{31,3}:\bq]=15$. 

For all other cases when \eqref{boundforK} gives a bound
$[\bk:\bq]\le t$ where $t>120$, we can similarly improve this bound using the Method A
applied to $\bff=\bff_{k,s}$. For example, for $k=23$ and $s=3$,
the inequality \eqref{boundforK} gives $N=[\bk:\bq]\le 3091$, and
\eqref{boundforKF} gives $[\bk:\bff_{23,3}]\le 281$. Applying the Method A,
we obtain $[\bk:\bff_{23,3}]\le 8$ and $[\bk:\bq]\le 11\cdot 8=88$. Surprisingly,
this strategy works in all bad cases.

Thus, our upper bound for degrees of fields from
$\F \Gamma^{(4)}_5(14)$ is $120$.

\medskip

Now, considering the graphs $\Gamma^{(4)}_4(14)$ again,
from \eqref{eq42N}, we similarly get
the inequalities
\begin{equation}
\frac{\ln{8}-\ln{\sin{\frac{\pi}{k}}}-
\ln{\sin{\frac{\pi}{r}}}}
{ \ln{2}-\frac{\ln{\gamma(k)}}{\varphi(k)}-\frac{\ln{\gamma(r)}}
{\varphi(r)}
}>\frac{\varphi([k,r])}{2\rho(k,r)}, \ \
\label{boundforkr}
\end{equation}

\begin{equation}
[\bk:\bff_{k,r}]\le \left[\frac{\ln{8}-\ln{\sin{\frac{\pi}{k}}}-
\ln{\sin{\frac{\pi}{r}}}}
{ \ln{2}-\frac{\ln{\gamma(k)}}{\varphi(k)}-\frac{\ln{\gamma(r)}}
{\varphi(r)}
}\left/\frac{\varphi([k,r])}{2\rho(k,r)} \right.\right],\,
\label{boundforKF4}
\end{equation}
and
\begin{equation}
N=[\bk:\bq]\le \left[\frac{\ln{8}-\ln{\sin{\frac{\pi}{k}}}-
\ln{\sin{\frac{\pi}{r}}}}
{ \ln{2}-\frac{\ln{\gamma(k)}}{\varphi(k)}-\frac{\ln{\gamma(r)}}
{\varphi(r)}
}\left/\frac{\varphi([k,r])}{2\rho(k,r)} \right.\right]\
\frac{\varphi([k,r])}{2\rho(k,r)}\ .
\label{boundforK4}
\end{equation}
for non-exceptional $k\ge r$ where $r=3,4,5$ and $k\ge 6$
(thus, one should replace $s$ by $r$ and $7$ by $8$ in our considerations
of $\Gamma^{(4)}_5(14)$ above).

Then exactly the same considerations as for $\Gamma^{(4)}_5(14)$
show that $[\bk:\bq]\le 120$ for all fields $\bk$ from $\Gamma_4^{(4)}(14)$
where $120$ is achieved for $k=31$ and $r=3$.

This finishes the proof of Theorem \ref{thdegrees}.

\section{A mirror symmetric finiteness conjecture\\
about reflective automorphic forms on\\
Hermitian symmetric domains of type IV}
\label{sec:refautforms}

In \cite{GrNik1}---\cite{GrNik4} some finiteness results and conjectures
about so called {\it reflective automorphic forms on symmetric domains
of type IV} (in classification by Cartan) were obtained and formulated.
They were considered as
mirror symmetric statements to finiteness results about
arithmetic hyperbolic reflection groups over $\bq$ and
corresponding reflective hyperbolic lattices over $\bz$.

Now finiteness results about arithmetic hyperbolic reflection groups
are established in full generality. Moreover, results of this paper can
be considered as some steps to classification in the future.
Respectively, it would be interesting to extend and formulate the
corresponding finiteness conjecture about reflective automorphic forms
on symmetric domains of type IV in full generality too. Let us do it.

Let $\bk$ be a totally real algebraic number field and $\bo$ its ring of
algebraic integers.

We recall that a lattice $L$ over $\bk$ is a finitely generated torsion-free
$\bo$-module $L$ with a symmetric bilinear form defined on $L$ with values in
$\bo$. Here $\bk$ is called the ground field of $L$, the number
$\dim L\otimes_\bo\bk$ is called the rank of $S$, and the absence of
torsion means that $S\subset L\otimes_\bo\bk$. We let $x\cdot y$
denote the value of the bilinear form on $L$ on the pair of elements
$x,y\in L$, and we let $x^2$ denote $x\cdot x$.

A lattice $S$ is said to be {\it hyperbolic} if the bilinear form
$S\otimes_\bo\br$ over $\br$ is indefinite for exactly one embedding
$\sigma^{(+)}:\bk\to \br$, and it is hyperbolic under this embedding, i. e.,
it has signature $(1,t_{(-)})$. For $\rk S\ge 3$ we let $\La(S)$ denote
the hyperbolic space of the dimension $\rk S -1$
which is canonically associated with the hyperbolic
lattice $S$ (and with the form $S\otimes_\bo \br$ under the embedding
$\sigma^{(+)}$):
$$
\La(S)=\{\br^+x\subset S\otimes_\bo \br\ |\ x^2>0\}_0\
$$
where $0$ means that we take a connected component.
The automorphism group $O^+(S)$ of the hyperbolic lattice $S$
is discrete and arithmetic in $\La(S)$, and it has fundamental domain
of finite volume.

A lattice $T$ is said to be {\it of IV type} if the bilinear form
$S\otimes_\bo\br$ over $\br$ is indefinite for exactly one embedding
$\sigma^{(+)}:\bk\to \br$, and
it has signature $(2,t_{(-)})$ for this embedding. For $\rk T\ge 3$ we let
$\Omega (T)$
denote the Hermitian symmetric domain of type IV (of dimension
$\rk T - 2$)
which is canonically associated with
the IV type lattice $T$ (and with the form $T\otimes_\bo \br$
under the embedding $\sigma^{(+)}$):
$$
\Omega (T)=\{\bc\omega\subset T\otimes \bc\ |\ \omega^2=0, \omega\cdot
\overline{\omega}>0\}_0\ .
$$
The automorphism group $O^+(T)$ of the IV type lattice $T$
is discrete and arithmetic in $\Omega(T)$, and it has fundamental domain of
finite volume.

For both types of lattices (hyperbolic or of IV type) $L$ we can define
reflections as follows.
Let $\delta \in L$, where $\sigma^{(+)}(\delta^2)<0$ and
$\delta^2\,|\,2(L\cdot \delta)$ (such elements
are called {\it roots} of $L$).
Then the formula
$$
s_\delta(x)=x-\frac{2(x\cdot \delta)}{\delta^2}\,\delta,\ x\in L,
$$
defines an involution $s_\delta$ of the lattice $L$ which is called
reflection relative to the root $\delta$ of $L$. We let $W(L)$ denote
the subgroup of $O(L)$ generated by all of the reflections of $L$.
Geometrically, in hyperbolic (respectively IV type) case the reflections
$s_\delta$ are precisely those automorphisms of $L$ which act as
reflections relative to hyperplanes of $\La(L)$ (respectively
to quadratic divisors
$D_\delta=\{\bc\omega\in \Omega(T)\ |\ \omega\cdot \delta=0\}$
of $\Omega(L)$). Obviously, they are orthogonal to the roots
$\delta$.

A hyperbolic lattice $S$ of the rank no less than three is called {\it
reflective} if $W(S)$ is a subgroup of $O(S)$ of finite index.

By Vinberg's arithmeticity criterion, any arithmetic hyperbolic
reflection group $W$ is a subgroup of finite index $W\subset W(S)$ for
one of reflective hyperbolic lattices $S$. The ground field $\bk$ of
$S$ then coincides with the ground field of $W$.

It was proved in \cite{Nik1} and \cite{Nik2} that for a fixed degree
$N=[\bk:\bq]$ of the ground field of $S$ and fixed rank $S\ge 3$
there exists only finitely many reflective hyperbolic lattices $S$
up to similarity (i.e. up to
multiplication of the form of $S$ by elements $k\in \bk$).
Since it is now established in full generality that
the $\rk S$ and degree $N$ are absolutely bounded, it follows that

\medskip

{\it there exist finitely many similarity classes of reflective hyperbolic
lattices}

\medskip

\noindent
(for all ranks $\ge 3$ and for all ground fields together).

Now let $T$ be a lattice of IV type. A holomorphic automorphic
form $\Phi$ of a positive weight on $\Omega(T)$ which is
automorphic relative to $O^+(T)$ is called {\it reflective
automorphic form of the lattice $T$} if the divisor of $\Phi$ is
union of quadratic divisors of $\Omega(T)$ which are orthogonal to
roots $\delta$ of reflections $s_\delta$ of $T$.  A lattice $T$ of
IV type is called {\it reflective} if it has at least one
reflective automorphic form $\Phi$.

\begin{conjecture} There exist finitely many similarity classes of
IV type reflective lattices of rank at least 5
(for all ranks $\ge 5$ and all ground fields together).
\end{conjecture}

For $\bk=\bq$ and respectively $\bo=\bz$ this conjecture
was formulated in \cite{GrNik1}--\cite{GrNik4}. Even in this case,
it seems, it is not established in full generality.

Some arithmetic hyperbolic reflection groups and some reflective automorphic
forms and corresponding hyperbolic and IV type reflective lattices
over $\bz$ are important in Borcherds proof \cite{Bor}
of Moonshine Conjecture by Conway and Norton \cite{CN}
which had been first discovered by John McKay.

We hope that similar objects over arbitrary number fields will find similar
astonishing applications in the future. At least, the results and conjectures
of this paper show that they are very exceptional even in this very
general setting.

\section{Appendix: Hyperbolic numbers (the review of \cite[Sec. 1]{Nik2})}
\label{sec:Appendix}

\subsection{Fekete's theorem}\label{subsec:fekete}
Here we review our results in \cite[Sec.1]{Nik2} and correct some
arithmetic mistakes (Theorems 1.1 and 1.2.2 in \cite{Nik2} which are similar
to Theorems \ref{thFekete2} and \ref{thhypnumberscond} here).
This mistakes are unessential for results of \cite{Nik2}.

The following
important theorem, to which this section is devoted, was obtained
by Fekete \cite{Fek}. Although Fekete considered (as we know) the
case of $\bq$ his method of proof can be immediately carried over
to totally real algebraic number fields.

\begin{theorem} (M. Fekete). Suppose that $\bff$ is a totally real
algebraic number field, and to every embedding $\sigma: \bff\to
\br$ there corresponds an interval $[a_\sigma,b_\sigma]$ in $\br$
and the real number $\lambda_\sigma>0$. Suppose that
$\prod_\sigma{\lambda_\sigma}=1$. Then for every nonnegative
integer $n$ there exists a nonzero polynomial $P_n(T)\in \bo[T]$
of degree no greater than $n$ over the ring of integers $\bo$ of
$\bff$ such that the following inequality holds for each $\sigma$:
\begin{equation}
\max_{x\in [a_\sigma,b_\sigma]}{|P_n^\sigma(x)|}\le \lambda_\sigma
 {|\discr \bff|^{1/(2[\bff:\bq])}
2^{n/(n+1)}(n+1)\left(\prod_\sigma {\frac{b_\sigma -
a_\sigma}{4}}\right)^{n/(2[\bff:\bq])}} \label{fek1}\ .
\end{equation}
\label{thFekete1}
\end{theorem}

\begin{proof} Suppose that $N=[\bff:\bq]$ and that
$\gamma_1,\dots\gamma_N$ is the basis for $\bo$ over $\bz$.
Suppose we are given a nonzero polynomial
$$
P_n(T)=\sum_{i=0}^{n}{\sum_{j=1}^{N}{\alpha_{ij}\gamma_jT^i}}\in
\bo[T]
$$
of degree no greater than $n$, where the $\alpha_{ij}\in \bz$
are not all zero. For every $\sigma:\bff\to \br$ we consider the
real functions $P_n^\sigma(x)$ on the interval
$[a_\sigma,b_\sigma]$.

We make the change of variables
$$
x=x(z)=\frac{b_\sigma+a_\sigma}{2}+\frac{b_\sigma-a_\sigma}{2}\cos{z}.
$$
If $z$ runs through $[0,\pi]$, then $x$ runs through
$[a_\sigma,b_\sigma]$. We also set
$Q^\sigma_n(z)=P_n^\sigma(x(z))$.

Since $Q_n^\sigma(z)$ is an even trigonometric polynomial, it
follows that
\begin{equation}
Q_n^\sigma(z)=\sum_{k=0}^{n}A_{k\sigma}\cos{kz},
\label{Qnsigma}
\end{equation}
where
$$
A_{k\sigma}=\frac{1}{\pi}\int_{-\pi}^{\pi}
{P_n^\sigma\left(\frac{b_\sigma+a_\sigma}{2}+\frac{b_\sigma-a_\sigma}{2}\cos{z}
\right)\cos{kz}\,dz}
$$
$$
=\sum_{i=0}^{n}{\sum_{j=1}^{N}\left(\frac{1}{\pi}\int_{-\pi}^{\pi}{\gamma_j^\sigma
\left(\frac{b_\sigma+a_\sigma}{2}+\frac{b_\sigma-a_\sigma}{2}\cos{z}
\right)^i\cos{kz}\,dz}\right)}\alpha_{ij}\,,
$$
if $k\ge 1$, and
$$
A_{0\sigma}=\frac{1}{2\pi}\int_{-\pi}^{\pi}
{P_n^\sigma\left(\frac{b_\sigma+a_\sigma}{2}+\frac{b_\sigma-a_\sigma}{2}\cos{z}
\right)\,dz}
$$
$$
=\sum_{i=0}^{n}{\sum_{j=1}^{N}\left(\frac{1}{2\pi}\int_{-\pi}^{\pi}
{\gamma_j^\sigma
\left(\frac{b_\sigma+a_\sigma}{2}+\frac{b_\sigma-a_\sigma}{2}\cos{z}
\right)^i \,dz}\right)}\alpha_{ij}\,.
$$
Thus,
\begin{equation}
A_{k\sigma}=\sum_{i=0}^{n}{\sum_{j=1}^{N}{c_{k\sigma
ij}\alpha_{ij}} \label{Aksigma}}
\end{equation}
are linear functions of the $\alpha_{ij}$, where
$$
c_{k\sigma ij}=\gamma_j^\sigma\cdot \frac{1}{\pi}\int_{-\pi}^{\pi}
{\left(\frac{b_\sigma+a_\sigma}{2}+\frac{b_\sigma-a_\sigma}{2}\cos{z}
\right)^i\cos{kz}\,dz}\, ,
$$
if $k\ge 1$, and
$$
c_{0\sigma ij}=\gamma_j^\sigma\cdot
\frac{1}{2\pi}\int_{-\pi}^{\pi}
{\left(\frac{b_\sigma+a_\sigma}{2}+\frac{b_\sigma-a_\sigma}{2}\cos{z}
\right)^i \,dz}\, .
$$
We note that, because of these formulas, $c_{k\sigma ij}=0$ for
$i<k$, and
$$
c_{k\sigma kj}=\gamma_j^\sigma\cdot
2\left(\frac{b_\sigma-a_\sigma}{4}\right)^k \,, \text{\ \  if
}k\ge 1
$$
($c_{0\sigma0j}=\gamma_j^\sigma$). Hence, if we order the indices
$k\sigma$ and $ij$ lexicographically, we find that the matrix of
the linear forms \eqref{Aksigma} is an upper block-triangular
matrix with the shown above $N\times N$ matrices $(c_{0\sigma 0 j})$ and
$(c_{k\sigma k j})$, $1\le k \le n$, on
the diagonal. It follows that its determinant is equal to
$$
\Delta=\det(\gamma_j^{\sigma})^{n+1}\cdot
2^{Nn}\left(\prod_\sigma{\frac{b_\sigma-a_\sigma}{4}}\right)^{n(n+1)/2}\,.
$$
Since $\prod_\sigma \lambda_\sigma^{n+1}=
\left(\prod_\sigma \lambda_\sigma  \right)^{n+1}=1$,
according to Minkowski's theorem on linear forms (see, for
example, \cite{BS}, \cite{Cas}), there exist $\alpha_{ij}\in \bz$, not
all zero, such that
$|A_{k\sigma}|\le
\lambda_\sigma |\Delta|^{1/N(n+1)}$,
and hence, by \eqref{Qnsigma},
$$
\max_{z}{\vert Q_n^\sigma(z)\vert }\le
\lambda_\sigma \cdot (n+1)|\Delta|^{1/N(n+1)}.
$$
Taking into account that $\det(\gamma_j^\sigma)^2=\discr \bff$, we
obtain the proof of the theorem.
\end{proof}

Taking $\lambda_\sigma=1$, we get a particular statement which we later use.

\begin{theorem} (M. Fekete). Suppose that $\bff$ is a totally real
algebraic number field, and to every embedding $\sigma: \bff\to
\br$ there corresponds an interval $[a_\sigma,b_\sigma]$ in $\br$.
Then for every nonnegative
integer $n$ there exists a nonzero polynomial $P_n(T)\in \bo[T]$
of degree no greater than $n$ over the ring of integers $\bo$ of
$\bff$ such that the following inequality holds for each $\sigma$:
\begin{equation}
\max_{x\in [a_\sigma,b_\sigma]}{|P_n^\sigma(x)|}\le
 {|\discr \bff|^{1/(2[\bff:\bq])}
2^{n/(n+1)}(n+1)\left(\prod_\sigma {\frac{b_\sigma -
a_\sigma}{4}}\right)^{n/(2[\bff:\bq])}} \label{fek1}\ .
\end{equation}
\label{thFekete2}
\end{theorem}

\subsection{Hyperbolic numbers}
\label{hypnumbers}
The totally real algebraic integers $\{\alpha\}$ which we consider
here are very similar to Pisot-Vijayaraghavan numbers \cite{Cas}, although
the later are not totally real.

\begin{theorem} Let $\bff$ be a totally real algebraic number field,
and let each imbedding $\sigma:\bff \to \br$ corresponds
to an interval $[a_\sigma,b_\sigma]$ in $\br$, where
$$
\prod_\sigma{\frac{b_\sigma-a_\sigma}{4}}<1.
$$
In addition, let the natural number $m$ and the intervals
$[s_1,t_1],\dots , [s_m,t_m]$ in $\br$
be fixed.

Then there exists a constant $N(s_i,t_i)$ such that, if $\alpha$
is a totally real algebraic integer
and if the following inequalities hold for the imbeddings
$\tau:\bff(\alpha)\to \br$:
$$
s_i\le \tau(\alpha)\le t_i,\ \text{for\ }\tau=\tau_1,\dots, \tau_m,
$$
$$
a_{\tau\vert \bff}\le \tau(\alpha)\le b_{\tau\vert \bff}\
\text{for\ }\tau\not=\tau_1,\dots\tau_m,
$$
then
$$
[\bff(\alpha):\bff]\le N(s_i,t_i).
$$
\label{thhypnumbers}
\end{theorem}

\begin{theorem}
Under the conditions of Theorem \ref{thhypnumbers}, $N(s_i,t_i)$
can be taken to be $N(s_i,t_i)=N(S)$, where $N(S)$ is the least natural
number solution of the inequality
\begin{equation}
nln(1/R)-M\ln(2n+2)-\ln{B}\ge \ln{S}.
\label{ineqforn}
\end{equation}
Here
\begin{equation}
M=[\bff:\bq],\ \ \
R=\sqrt{\prod_\sigma{\frac{b_\sigma-a_\sigma}{4}}},
\label{ineqforMR}
\end{equation}
\begin{equation}
B=\sqrt{\vert \discr \bff\vert},\ \ \
S=\prod_{i=1}^{m}{\left( 2e r_i(b_{\sigma_i}-a_{\sigma_i})^{-1}\right)},
\label{ineqforBS}
\end{equation}
where
$\sigma_i=\tau_i\vert\bff$ and $r_i=\max{\{\vert t_i-a_{\sigma_i}\vert,\
\vert b_{\sigma_i}-s_i\vert\}}$.

Asymptotically,
$$
N(s_i, t_i)\sim \frac{\ln S }{\ln {(1/R)}}.
$$
\label{thhypnumberscond}
\end{theorem}

\begin{proof} We use the following statement.

\begin{lemma} Suppose that $Q_n(T)\in \br[T]$ is a non-zero polynomial
over $\br$ of degree no greater than $n>0$,
$a<b$ and $M_0=\max_{[a,b]}{\{\vert Q_n(x)\vert\}}$. Then for $x\ge b$
$$
|Q_n(x)|\le \frac{M_0(x-a)^n n^n }{\left((b-a)/2\right)^n n!}<
$$
$$
\frac{M_0(x-a)^n e^n}{((b-a)/2)^n\sqrt{2\pi n}} <
\frac{M_0(x-a)^n e^n}{((b-a)/2)^n} .
$$
\label{lemLagr}
\end{lemma}

\begin{proof} Let $\alpha_0<\alpha_1<\cdots <\alpha_n$.
Then we have the Lagrange interpolation formula
$$
Q_n(x)=\sum_{i=0}^{n}{Q_n(\alpha_i)F_i(x)}
$$
where
$$
F_i(x)=\frac{(x-\alpha_0)(x-\alpha_1)\cdots
(x-\alpha_{i-1})(x-\alpha_{i+1})\cdots (x-\alpha_n)}
{(\alpha_i-\alpha_0)(\alpha_i-\alpha_1)\cdots
(\alpha_i-\alpha_{i-1})(\alpha_i-\alpha_{i+1})\cdots(\alpha_i-\alpha_n)}.
$$
Taking $\alpha_i=a+i(b-a)/n$, $0\le i\le n$, we obtain for $x\ge b$ that
$$
|Q_n(x)|\le \frac{M_0(x-a)^n}{((b-a)/n)^n}\sum_{i=0}^{n}\frac{1}{i!(n-i)!}=
\frac{M_0(x-a)^n2^n}{((b-a)/n)^n n!}.
$$
By Stirling formula, $n!=\sqrt{2\pi n}(n/e)^ne^{\lambda_n}$ where
$0<\lambda_n<1/(12n)$. Thus, $n^n/n!< e^n/\sqrt{2\pi n}<e^n$. It
follows the statement.
\end{proof}

We continue the proof of theorems.

For given $n$ we consider the polynomial $P_n(T)\in \bo[T]$ whose
existence is ensured by Fekete's theorem \ref{thFekete2}. Setting
$N=[\bff(\alpha):\bff]$ and $M=[\bff:\bq]$, we use Fekete's
theorem and the lemma to conclude that
$$
|\prod_\tau{\tau(P_n(\alpha))|}=\prod_{\tau}{|P_n^\tau(\tau(\alpha))|}=
\prod_{\tau\not=\tau_i}{P_n^\tau(\tau(\alpha))}
\prod_{i=1}^{m}{|P_n^{\tau_i}(\tau_i(\alpha))|}
$$
$$
\le \prod_{\tau\not=\tau_i}{\max_{[a_{\tau|\bff},b_{\tau|\bff}]}{|P_n^\tau(x)|}}
\prod_{i=1}^m{\max_{[s_i,t_i]}{|P_n^{\sigma_i}(x)|}}
$$
$$
\le\left(|\discr \bff|^{1/(2M)}\cdot 2\cdot (n+1)R^{n/M}\right)^{NM}
\prod_{i=1}^{m}{\frac{r_i^ne^n}{((b_{\sigma_i}-a_{\sigma_i})/2)^n}}
$$
$$
=R^{nN}B^N\cdot S^n\cdot (2n+2)^{MN}.
$$
Since $R<1$, there exists $n_0$ large enough so that
\begin{equation}
R^{n_0}\cdot B\cdot(2n_0+2)^{M}\le \frac{1}{S}.
\label{eqn0}
\end{equation}
Then if $N>n_0$, we find that
$$
R^{n_0N}\cdot B^N\cdot S^{n_0}(2n_0+2)^{MN}\le S^{n_0-N}<1,
$$
since $S>1$. From this and the above chain of inequalities we have
$$
|\prod_\tau{\tau(P_{n_0}(\alpha))}|<1.
$$
But
$$
\prod_\tau{\tau(P_{n_0}(\alpha))}=N_{\bff(\alpha)/\bq}(P_{n_0}(\alpha))\in \bz,
$$
and hence $P_{n_0}(\alpha)=0$. Consequently, $N\le n_0$, and we have
obtained a contradiction.
We have thereby proved that $N\le n_0$, where $n_0$ is a natural
number solution of \eqref{eqn0}.
The inequality \eqref{eqn0} is obviously equivalent to
$$
n_0\ln{(1/R)}-M\ln{(2n_0+2)}-\ln{B}\ge \ln{S},
$$
and this completes the proof of the theorems.
\end{proof}


V.V. Nikulin \par Deptm. of Pure Mathem. The University of
Liverpool, Liverpool\par L69 3BX, UK; \vskip1pt Steklov
Mathematical Institute,\par ul. Gubkina 8, Moscow 117966, GSP-1,
Russia

vnikulin@liv.ac.uk \ \ vvnikulin@list.ru

\end{document}